\documentstyle[12pt]{article}
\renewcommand{\theequation}{\thesection.\arabic{equation}}

\begin{document}
\renewcommand{\thefootnote}{\fnsymbol{footnote}}
\newpage
\pagestyle{empty}
\setcounter{page}{0}

\vfill
\begin{center}

{\LARGE {\bf {\sf 
Canonical factorization and diagonalization of Baxterized braid matrices: Explicit
constructions and applications. }}} \\[0.8cm]
{\large A.Chakrabarti

{\em 

Centre de Physique Th\'eorique\footnote{Laboratoire Propre 
du CNRS UPR A.0014}, Ecole Polytechnique, 91128 Palaiseau Cedex, France.\\
e-mail chakra@cpht.polytechnique.fr}}

\end{center}

\smallskip

\smallskip 

\smallskip

\smallskip

\smallskip

\smallskip 

\begin{abstract}
Braid matrices $\hat{R}(\theta)$, corresponding to vector representations, are spectrally
decomposed obtaining a ratio $f_{i}(\theta)/f_{i}(-\theta)$ for the coefficient of each
projector
$P_{i}$ appearing in the decomposition. This directly yields a factorization
$(F(-\theta))^{-1}F(\theta)$ for the braid matrix, implying also the relation
$\hat{R}(-\theta)\hat{R}(\theta)=I$.This is achieved for
$GL_{q}(n),SO_{q}(2n+1),SO_{q}(2n),Sp_{q}(2n)$ for all $n$ and also for various other
interesting cases including the $8$-vertex matrix.We explain how the limits $\theta
\rightarrow \pm \infty$ can be interpreted to provide factorizations of the standard
(non-Baxterized) braid matrices. A systematic approach to diagonalization of projectors
and hence of braid matrices is presented with explicit constructions for
$GL_{q}(2),GL_{q}(3),SO_{q}(3),SO_{q}(4),Sp_{q}(4)$ and various other cases such as the
$8$-vertex one. For a specific nested sequence of projectors diagonalization is
obtained for all dimensions. In each factor $F(\theta)$ our diagonalization again
factors out all dependence on the spectral parameter $\theta$ as a diagonal matrix. The
canonical property implemented in the diagonalizers is mutual orthogonality of the
rows. Applications of our formalism to the construction of $L-$operators and transfer
matrices are indicated. In an Appendix our type of factorization is compared to
another one proposed by other authors.  
\end{abstract}

\vfill
\newpage

\pagestyle{plain}

\section{Introduction}
 Let $\hat{R}(\theta)$ be a braid matrix Baxterized with a spectral parameter $\theta$
and satisfying, in standard notations, 
\begin{equation}
\hat{R}_{12}(\theta)\hat{R}_{23}(\theta+\acute{\theta)}\hat{R}_{12}(\acute{\theta)}
=\hat{R}_{23}(\acute{\theta})\hat{R}_{12}(\theta +\acute{\theta})\hat{R}_{23}(\theta)
\end{equation}

Here, apart from $\theta$, $\hat{R}(\theta)$ can depend on other parameters such as $q$,
which will not always be denoted explicitly. Vector representations with $N^{2}\times
N^{2}$ braid matrices are implied in all cases. The corresponding $YB$ (Yang-Baxter)
martix is
\begin{equation}
{R}(\theta)=P\hat{R}(\theta)
\end{equation}
where the permutation matrix $P$ is defined to be (with $i=(1,2,...N)$)
 \begin{equation}
P=\sum_{ij} E_{ij}\bigotimes E_{ji}
\end{equation}
The matrix $E_{ij}$ has zero elements except for a single unit one at $(ij)$.

 We assume that the polynomial equation of minimal degree satisfied by $\hat{R}(\theta)$
has distinct roots. When this holds $\hat{R}(\theta)$ can be spectrally decomposed on a
basis of projectors $P_i$, satisfying 

\begin{equation}
P_{i}P_{j}=\delta_{ij}, \qquad  \sum_{i} P_i = I_{N^2\times N^2}
\end{equation}
 Suppressing arguments for the time being, if (with $k_i\neq k_j$ for $i\neq j$)
\begin{equation}
\prod_{i=1}^p (\hat R -k_{i}I)=0
\end{equation}

then defining
\begin{equation}
P_i= \prod_{j\neq i}\frac{(\hat R -k_{j}I)}{(k_i -k_j)}, \qquad (i=1,2,...p)
\end{equation}

the set $P_i$ can be shown to satisfy $(1.4)$ and one obtains

\begin{equation}
\hat R=\sum_{i}^p k_i P_i
\end{equation}

On the other hand,given $(1.7)$ one obtains $(1.5)$. The $P$'s on the right can, in
general, depend on parameters such as $q$. {\it But in all cases they will be independent
of the spectral parameter} $\theta$. In $R(\theta)$ all $\theta$-dependence is to be
found in the coefficients $k_i$. This is consistent with $(1.6)$ and $(1.7)$ and is
fundamental for the considerations below.

 In all cases to be considered, not only we will obtain explicit spectral decomposition
of $\hat{R}(\theta)$, but also {\it a specific factorized form} of each $k_i$

\begin{equation}
k_i(\theta)=\frac{f_i(\theta)}{f_i(-\theta)}
\end{equation}
when
\begin{equation}
\hat{R}(\theta) =\sum_{i} k_i(\theta)P_i=\sum_{i}\frac{f_i(\theta)}{f_i(-\theta)}P_i
\end{equation}

 This will be our first major step.

  The number of projectors and their matrix elements are specific to the case considered.
But they always satisfy $(1.4)$. In $[1]$, $(1.9)$ has been obtained explicitly for
$GL_q(n), SO_q(2n+1),SO_q(2n)$ and $Sp_q(2n)$ for all $n$.The results are recapitulated
in $Secs.2,3$.In $Secs.4,5,6,7$ we obtain $(1.9)$ for various interesting cases,
including the $8-$vertex matrix.

  An evident, but for us crucial, consequence of $(1.4)$ is that for well-defined and
mutually commuting but otherwise {\it arbitray} coefficients $(a_i,b_i)$,
\begin{equation}
(\sum_{i} a_i P_i)(\sum_{i} b_i P_i)=(\sum_{i} a_ib_i P_i)=(\sum_{i} b_i P_i)(\sum_{i}
a_i P_i)
\end{equation}

Hence, once a spectral decomposition $(1.7)$ has been obtained $\hat R$ can be expressed
as a product of arbitrary number of factors
\begin{equation}
\hat R=\prod_{n}(\sum_{i} k_i^{(n)} P_i), \qquad  (\prod_{n} k_i^{(n)} =k_i)
\end{equation}

 Of particular interest to us is the factorization
$$\hat{R}(\theta)
=\sum_{i}\frac{f_i(\theta)}{f_i(-\theta)}P_i=(\sum_{i}f_i^{-1}(-\theta)P_i)(\sum_{i}f_i(\theta)P_i)$$

\begin{equation}
=(F(-\theta))^{-1}F(\theta)
\end{equation}
where
\begin{equation}
F(\theta)= \sum_{i}f_i(\theta)P_i
\end{equation}

This implies the so-called "unitarity"

\begin{equation}
\hat R(-\theta)\hat R(\theta) = I_{N^2\times N^2}
\end{equation}

  One obtains from $(1.2)$, ,since $P^2=I$,

\begin{equation}
R(\theta)=(P(F(-\theta))^{-1}P)PF(\theta)=(F_{21}(-\theta))^{-1}P F_{12}(\theta)
\end{equation}

 In $(1.12)$ and $(1.15)$ the key feature is the change of sign of $\theta$ in $F^{-1}$.

  Other interesting choices are possible. Thus, for example, defining
\begin{equation}
\acute {F}(\theta)=
\sum_{i}\Biggl(\frac{f_i(\theta)}{f_i(-\theta)}\Biggr)^{\frac{1}{2}}P_i=(\acute
{F}(-\theta))^{-1}
\end{equation}

 one obtains 
\begin{equation}
\hat R(\theta)=(\acute {F}(-\theta))^{-1}\acute{F}(\theta) =(\acute{F}(\theta))^2
\end{equation}

and
\begin{equation}
 R(\theta)=(\acute {F}_{21}(-\theta))^{-1}P\acute{F}_{12}(\theta)
=(\acute{F}_{21}(\theta))P\acute{F}_{12}(\theta)
\end{equation}

  Here, even for real $\hat R(\theta)$, for certain domains of $\theta$ the factor 
$\acute{F}(\theta)$ can be complex.

  Compare $(1.15)$ and $(1.18)$ to a Drinfeld twist $[2,3,4]$ of $P$

\begin{equation}
 R''(\theta)=( F''_{21}(\theta))^{-1}PF''_{12}(\theta)
\end{equation}

 In $(1.15)$ there is $(-\theta)$ on the left and in $(1.18)$ there is no inversion of
$\acute F_{21}(\theta)$. $P$ satisfies the $YB$ equation with the trivial $\hat R=P^2=I$
for the braid matrix. ( $P$ also satisfies the braid equation with $ R=P^2=I$.) The
properties of $ R''(\theta)$ will depend on those of $F''(\theta)$ ( such as cocycle
conditions). In our case, since one {\it starts} from solutions of $(1.1)$ one does not
have to verify if $F(\theta)$ and $\acute {F}(\theta)$ satisfy suitable constraints, so
far as the braid equation is concerned.

  The present situation may also be compared to "contraction" of $YB$ matrices to
non-standard, Jordanian forms. Without even trying to explain the terminology we refer
to two $[5,6]$ of our series of relevant papers (where original sources are cited). We
mention this only to point out that , as compared to $(1.15),(1.18),(1.19)$ the role of
$P$ is, so to say, reversed. For the non-standard case $R$ is a Drinfeld twist of $I$,

\begin{equation}
R=(F_{21})^{-1}F,  \qquad \hat R = F^{-1}PF.
\end{equation} 

 The nontrivial matrix $\hat R$ is now "triangular" since from $(1.20)$
\begin{equation}
 (\hat R)^2 = I
\end{equation} 

  The ambiguities arising in factorizing (compare $(1.13)$ and $(1.16)$), or in defining
$f_i(\theta)$ for a given $k_i(\theta)$ in $(1.8)$, become particularly relevant in
considering the limits $\theta \rightarrow \pm \infty$. From $(1.9)$ one has evidently

\begin{equation}
\hat R(0) = \sum_{i}P_i = I
\end{equation}

 It will be seen in the following sections that in each case for $\theta \rightarrow \pm
\infty$ one obtains the standard ( non-Baxterized) braid matrices ($\hat R$ and the
inverse) satisfying
\begin{equation}
\hat{R}_{12}\hat{R}_{23}\hat{R}_{12}
=\hat{R}_{23}\hat{R}_{12}\hat{R}_{23}
\end{equation}

This equation can be considered, consistently with $(1.1)$ as the limiting form (with
arguments suppressed) as both $(\theta,\acute \theta) \rightarrow +\infty$, say. If one
denotes 
\begin{equation}
lim_{\theta \rightarrow + \infty} \hat R(\theta)=\hat R
\end{equation}

then consistently with $(1.14)$
\begin{equation}
lim_{\theta \rightarrow - \infty} \hat R(\theta)={\hat R}^{-1}
\end{equation}

 {\it In these limits special features arise concerning factorizations.} It is helpful to
consider a simple but frequently encountered example. Suppose that for some
$f_i(\theta)$ one has (dropping the index $i$ and setting $q=e^h$)

\begin{equation}
\frac {f(\theta)}{f(-\theta)}=\frac {sinh(h-\theta)}{sinh(h+\theta)}
\end{equation}
 The evident singularity at $\theta=-h$ can be excluded by definition from the domain of
$\theta$. Now, as $\theta \rightarrow \pm \infty$
\begin{equation}
\frac {f(\theta)}{f(-\theta)} \rightarrow -q^{\mp 2}
\end{equation}

 But what about the factor $f(\theta)$ ? How does it behave when separated in the factor
$F(\theta)$ or $\acute F(\theta)$ exhibited before ?

  (1): For the choice
 \begin{equation}
f(\pm\theta) = sinh(h\mp\theta)
 \end{equation}
separately both $f(\theta)$ and $f(-\theta)$ both diverge.

  (2): The choice
\begin{equation}
f(\theta) = \Biggl(\frac{sinh(h-\theta)}{sinh(h+\theta)}\Biggr)^\frac{1}{2}
           = \frac{sinh(h-\theta)}{{(sinh(h-\theta)sinh(h+\theta))}^{\frac{1}{2}}}
 \end{equation}

gives consistently with $(1.27)$, finite but imaginary limits

\begin{equation}
f(\theta) \rightarrow \pm iq^{\mp 1}, \qquad f(-\theta) \rightarrow \mp iq^{\pm 1}
\end{equation}

(3): But more generally than in $(1.29)$ one may choose for any well-defined $y(\theta)$,

\begin{equation}
f(\theta) = \frac{sinh(h-\theta)}{{(y(\theta)y(-\theta))}^{\frac{1}{2}}}
 \end{equation}

Setting, for example,

\begin{equation}
f(\theta)  = \frac{sinh(h-\theta)}{{(cosh(h-\theta)cosh(h+\theta))}^{\frac{1}{2}}}
 \end{equation}

  one obtains real, finite limits consistent with $(1.28)$
\begin{equation}
\theta \rightarrow +\infty, \qquad f(\pm \theta) \rightarrow \mp q^{\mp1}
\end{equation}
 with an evident analogous result for $\theta \rightarrow -\infty$.

      We will assume that each $f_{i}(\theta)$ in $F(\theta)$ has thus been suitably
defined ( choosing an appropriate $y(\theta)$ ). Then even $\hat R$ satisfying $(1.23)$
can be considered to be factorized as in $(1.12)$, the implicit spectral parameter not
being exhibited in the limits $\theta \rightarrow \pm \infty$. In this sense, the
unitarity $(1.14)$ can still be considered to be implicit. Note that even if each
$f_{i}(\theta)$ has limits analogous to $(1.33)$ with different powers of $q$, one ${\it 
cannot}$ express the factorization as $(F(-q))^{-1}F(q))$ since the projectors,in    
 general, are $q-$dependent ( though always independent of $\theta$ ). It is essential
to think in terms of $\theta$ even when it is , in the limits above, implicit. The
implementation of a spectral parameter, the passage from $(1.23)$ to $(1.1)$, renders
many aspects more complex. But it also provides an extra margin of maneuvre, making
possible the canonical factorization $(1.12)$ possible whose interest will be studied
later on.

   Our factorizations are directly based on the resolution $(1.9)$. Other classes of
factorizations can also be envisaged. One such class with upper and lower triangular
factors for the $YB$ matrix $R(\theta)$, leading to interesting properties has been
studied by Maillet et al. in a series of papers $[7,8,9]$. This formalism is compared
with ours in $App.A$.

    Since all braid matrices studied $(Secs.2,3,4,5,6,7)$ are systematically found to
lead to spectral decompositions with each coefficient of the form  $(1.8)$ and $(1.9)$ a
more general study of such forms should be of interest. Here we will limit our
observations to the following feature. Let
\begin{equation}
\hat{R}'(\theta)=\sum_{i}\frac{g_i(\theta)}{g_i(-\theta)}P_i
\end{equation}
where, apart from being well-defined, the $g$'s are as yet arbitrary.In general,
$\hat{R}'(\theta)$ does $\it not$ satisfy $(1.1)$. But defining
\begin{equation}
\hat{H}(\theta)=\sum_{i}\frac{g_i(\theta)}{f_i(\theta)}P_i \equiv \sum_{i}h_i(\theta)P_i
\end{equation}
where $f_{i}(\theta)$ corresponds to $(1.9)$, a solution  $\hat{R}(\theta)$ of $(1.1)$
on the ${\it same}$ basis of projectors
\begin{equation}
\hat{R}'(\theta)={(H(-\theta))}^{-1}\hat{R}(\theta) H(\theta)
\end{equation}

  Any two matrices decomposable on the same spectral basis ( satisfying $(1.4)$ ) are
always related as above.

      Substituing in $(1.1)$ 
\begin{equation}
\hat{R}(\theta)=H(-\theta){\hat{R}}'(\theta) {H(\theta)}^{-1}
\end{equation}
one can rephrase $(1.1)$ in terms of $\hat {R}'$ and $H$. One obtains
\begin{equation}
\hat{R}'_{12}(\theta)X_{1}\hat{R}'_{23}(\theta+\acute{\theta)}X_{2}\hat{R}'_{12}(\acute{\theta)}X_{3}
=X_{4}\hat{R}'_{23}(\acute{\theta})X_{5}\hat{R}'_{12}(\theta
+\acute{\theta})X_{6}\hat{R}'_{23}(\theta)
\end{equation}
where 
$$ X_{1} = (H_{12}(\theta))^{-1}H_{23}(-\theta -\acute{\theta}) $$

and so on.

 Now, along with the properties of $f_i(\theta)$, those of  $g_i(\theta)$ will determine
the content of this equation for $\hat {R}'$. Further study in this direction is beyond
the scope of this work.

      Our first basic step is the systematic expression of the solutions of $(1.1)$ in
the form $(1.9)$. The next major one is the simultaneous diagonalization of each
projector $P_i$ in $(1.9)$ and hence of $\hat {R}(\theta)$. Our approach is presented
step by step in $Sec.9$. Explicit examples of diagonalizations of lower dimensional
cases of $Sec.2$ and $Sec.3$ $(Gl_q (2),Gl_q (3),SO_q(3),SO_q(4),Sp_q(4))$ are
collected together in $App.B$. At the end of $Secs.(4,5,6,7)$ the diagonalizations are
presented explicitly for each case. Our $Sec.8$ is an exception, where a nested
sequence of projectors with simple, attactive features is presented for arbitrary
dimensions without constructing explicit solutions of the braid equation. On the
contrary, here the diagonalizer is obtained quite simply for arbitrary dimensions.

  A {\it canonical} feature sought for in our formalism is the {\it mutual
orthogonality of the rows } of the matrix diagonalizing $\hat R(\theta)$. The elegant
and useful consequences of such a feature are pointed out. In the factorized form, our
diagonalization factors out again in each factor all $\theta$-dependence as a diagonal
matrix. 

   Applications of our spectral decompositions and diagonalizations to the construction
of $L$-operators and to transfer matrices are discussed respectively in $Sec.10$
and $Sec.11$.

\section{Factorization of braid matrices of $GL_q(N), SO_q(N)$ and $Sp_q (N)$ :}
\setcounter{equation}{0}
\renewcommand{\theequation}{2.\arabic{equation}}
 We recapitulate below the relevant essential results of $[1]$. The standard
$q-$dependent $N^2 \times N^2$ projectors $[10]$ are assumed to be known. For $Sp_q$
always $N=2n$.

  Same notations will be used for projectors in different cases though they are
different. The overall normalizing factor for $\hat R(\theta)$ is chosen to obtain $1$
for the element $(11)$ at top left. (See however $Sec.7$.)

  For $GL_q (N)$ one has two projectors $(P_+,P_-)$ satisfying $(1.4)$. For $\hat
R(\theta)$ satisfying $(1.1)$, setting $h=ln q$, one obtains 

$$ \hat R(\theta) = P_+ +\frac {sinh(h-\theta)}{sinh(h+\theta)} P_-$$

$$=(P_+ +(sinh(h+\theta))^{-1} P_-)(P_+ +{sinh(h-\theta)} P_-)$$
\begin{equation}
\equiv (F(-\theta))^{-1}F(\theta)
\end{equation}

 To illustrate $(1.12)$ we have implemented one simple possible choice for $F(\theta)$.
Ambiguities discussed in $Sec.1$ (from $(1.28)$ to $(1.32)$) {\it are always present in
this and other examples to follow}. This statement will not be repeated in successive
sections.

   For $SO_q(N)$, for $N=(2n+1)$ and also for $N=2n$, one has a basis of three
projectors $(P_+,P_-,P_0)$ and two possibilities:

\begin{equation}
 \hat R(\theta) = P_+ +\frac {sinh(h-\theta)}{sinh(h+\theta)} P_-                  
        +\frac {cosh(\frac{N}{2}h-\theta)}{cosh(\frac{N}{2}h+\theta)} P_0
\end{equation}

or
\begin{equation}
 \hat R(\theta) = P_+ +\frac {sinh(h-\theta)}{sinh(h+\theta)} P_-                  
        +\frac {sinh((\frac{N}{2}-1)h-\theta)
sinh(h-\theta)}{sinh((\frac{N}{2}-1)h+\theta) sinh(h+\theta)} P_0
\end{equation}

 For $Sp_q (2n)$ one obtains
 
\begin{equation}
 \hat R(\theta) = P_+ +\frac {sinh(h-\theta)}{sinh(h+\theta)} P_-                  
        +\frac {sinh((n+1)h-\theta)}{sinh((n+1)h+\theta)} P_0
\end{equation}

or
\begin{equation}
 \hat R(\theta) = P_+ +\frac {sinh(h-\theta)}{sinh(h+\theta)} P_-                  
 +\frac {cosh(nh-\theta)sinh(h-\theta)}{cosh(nh+\theta) sinh(h+\theta)} P_0
\end{equation}

  The expressions for $F(\theta)$ and $(F(-\theta))^{-1}$ are evident in each case. See
however the remarks below $(2.1)$.
 For each case
\begin{equation}
\hat R(0) =I
\end{equation}

 For $\theta \rightarrow \pm\infty$, carefully taking limits, one respectively obtains:

 For $GL_q(N)$

\begin{equation}
\hat R = P_+ - q^{\mp 2}P_-
\end{equation}
 
 For $SO_q(N)$

\begin{equation}
\hat R = P_+ - q^{\mp 2}P_- +q^{\mp N}P_0
\end{equation}

 For $Sp_q(N)$

\begin{equation}
\hat R = P_+ - q^{\mp 2}P_- - q^{\mp (N+2)}P_0
\end{equation}

  These are the standard ( non-Baxterized ) braid matrices $[10]$ satisfying $(1.23)$.
Concerning factorization see the relevant discussion in $Sec.1$ ( from $(1.23)$ to
$(1.33)$ ).

 At the end of $Sec.4$ of $[1]$ it has been pointed out that for $q=1 (h=0)$ all these
matrices reduduce to ones with constant elements satisfying

\begin{equation}
{\hat R}^2 =I
\end{equation}
 
They amount to twists of $I$ with constant matrices. This situation is to be contrasted
with the corresponding one in $Sec.3$.

\section{A new class of braid matrices for $SO_q(N)$ and $Sp_q(N)$:}
\setcounter{equation}{0}
\renewcommand{\theequation}{3.\arabic{equation}}

 This was presented in $Sec.4$ of $[1]$. The solution for $SO_q(3)$ appeared already in
$[11]$. The structure $(1.9)$ is again present and hence also the factorization $(1.12)$.

  We recapitulate:

  Define
\begin{equation}
d= (1+\epsilon [N-\epsilon])^{-1} \qquad = \biggl (1+\epsilon \frac 
{q^{N-\epsilon}-q^{-N+\epsilon}}{q-q^{-1}} \biggr )^{-1}
\end{equation}

 where for $SO_q(N)$

         $$\epsilon =1, \qquad N=3,4,...$$
 
and for $Sp_q(N)$

     $$\epsilon = -1, \qquad N=4,6,...$$

  Define also

\begin{equation}
tanh \eta = \sqrt {1-4d^2}
\end{equation}

 The reality of $\eta$ is implied by $(3.1)$ since $4d^2 <1$.
  Now with such an $\eta$, 
\begin{equation}
\hat R(\theta) = P_{+} +P_{-} + \frac {sinh(\eta -\theta)}{sinh(\eta +\theta)}P_{0}
\end{equation}
\begin{equation}
 = I+\biggl(\frac{sinh(\eta-\theta)}{sinh(\eta+\theta)}- 1
\biggr)P_0                          
\end{equation}

can be shown $[1]$ to satisfy $(1.1)$. The promised srtucture is explicit in 
$(3.3)$.

     One has as usual

\begin{equation}
\hat R(0) = I                       
\end{equation}

and for $\theta \rightarrow \pm \infty$ 

\begin{equation}
\hat {R}(\pm \infty) =P_+ + P_- -e ^{\mp 2\eta}P_0                          
\end{equation}

\begin{equation}
 =I -(1+e ^{\mp 2\eta})P_0                          
\end{equation}
where
\begin{equation}
 e ^{\mp 2\eta} = \frac {1\mp \sqrt {1-4d^2}}{1\pm \sqrt{1-4d^2}}                        
\end{equation}

 These provide a new class of ( non-Baxterized ) braid matrices ${\hat R}^{\pm1}$
satisfying $(1.23)$. The $R^{\pm1}=P{\hat R}^{\pm1}$ are new solutions of the $YB$
equation for $SO_q$ and $Sp_q$ with $\epsilon$ and $N$ as given below $(3.1)$.

 Moreover, from $(3.1)$, for $q=1$

\begin{equation}
d= \frac {\epsilon}{N}
\end{equation}

and from $(3.2)$, for $q=1$,
\begin{equation}
(tanh \eta)_{(q=1)} = \sqrt {1-\frac{4}{N^2}} \qquad \equiv tanh \hat{\eta}
\end{equation}
 Denoting 
\begin{equation}
(P_0)_{(q=1)} = \hat P_0
\end{equation}
 and so on, one obtains from $(3.3)$ and $(3.7)$ respectively

\begin{equation}
(\hat R(\theta))_{(q=1)} =\hat{P}_+ +\hat{P}_-
+\frac{sinh(\hat{\eta}-\theta)}{sinh(\hat{\eta}+\theta)}\hat{P}_0                          
\end{equation}
and
\begin{equation}
(\hat R^{(\pm1)})_{(q=1)} =\hat{P}_+ +\hat{P}_- - e
^{\mp2\hat{\eta}}\hat{P}_0                          
\end{equation}

 The braid matrix $(3.13)$ satisfies a nontrivial Hecke condition 

\begin{equation}
\bigl(\hat R - I \bigr)\bigl(\hat R + e^{-2\hat{\eta}} I \bigr) =0
\end{equation}

and cannot be twisted back to $I$. This situation should be compared to $(2.10)$  and
the remarks that follow $(2.10)$.

  Note that we are {\it not} expanding in powers of $h(=ln q)$ to extract the so-called
"classical" $r$-matrix. We are directly setting $q=1$ and yet getting quite nontrivial
results.

\section{Two exotic cases $(S03,S14)$:}
\setcounter{equation}{0}
\renewcommand{\theequation}{4.\arabic{equation}}

  Two special braid matrices arising in the classification of $4\times 4$ $YB$ matrices
of $[12]$ were Baxterized in $[13]$. Other aspects were already studied in previous
papers of the series $[14]$. Some "exotic" features are briefly recapitulated below in
the present context:

 $\bullet$  Complex projectors for $S03$ ( for real $\hat R$ )

 $\bullet$  Extended freedom of parametrization for $S14$. 

 Our solutions presented in $Sec.3$ can be considered to be an exotic class in arbitrary
dimensions $(N^2 \times N^2, N \geq 3)$. For even $N$ one has two types,  exotic
orthogonal and exotic symplectic.

\smallskip

   $ S03:$

\smallskip

              The braid matrix 

\begin{equation}
\hat R =  \pmatrix{
   1 &0 &0 &1 \cr
   0 &1 &-1 &0 \cr
   0 &1 &1 &0 \cr
   -1 &0 &0 &1                                                   
              }.
\end{equation}

  satisfies 

\begin{equation}
(\hat R -(1+i)I)(\hat R -(1-i)I) =0
\end{equation}

 The corresponding projectors 

\begin{equation}
P_{(\pm)}=\frac{1}{2}(I \pm i(\hat R -I))
\end{equation}
 provide the spectral decomposition

\begin{equation}
\hat R = (1-i)P_{(+)} + (1+i)P_{(-)}
\end{equation}

  Altering suitably the normalization of $[13]$ gives thne Baxterization (with
$z=e^\theta$)
\begin{equation}
\hat{R}(z) = \biggl(\frac{f(z)}{f(z^{-1})}\biggr)^{\frac{1}{2}}P_{(+)} +
\biggl(\frac{f(z^{-1})}{f(z)}\biggr)^{\frac{1}{2}} P_{(-)}
\end{equation}
where 
\begin{equation}
f(z)= (z+z^{-1}) +i(z-z^{-1})
\end{equation}

 Thus we obtain the form$(1.9)$ and $(1.12)$ follows.

 One can rewrite $(4.5)$ in the explicitly real form
\begin{equation}
\hat R (z) = (z ^2 +z^{-2})^{-\frac{1}{2}} \bigl({(\sqrt{2}z)^{-1}}\hat R
+\sqrt{2}z{\hat R}^{-1}\bigr)
\end{equation}
and verify again

\begin{equation}
\hat {R}(z^{-1})\hat {R}(z) = I
\end{equation}

The unitary matrix $M$, where

\begin{equation}
\sqrt{2}M=  \pmatrix{
   1 &0 &0 &i\cr
   0 &1 &-i &0 \cr
   0 &-i &1 &0 \cr
   i &0 &0 &1                                                   
              }.
\end{equation}

diagonalizes $P_{(\pm)}$ giving
\begin{equation}
M\hat R M^{-1}= diag(1-i,1-i,1+i,1+i)
\end{equation}

\begin{equation}
(z^2 +z^{-2})^{\frac{1}{2}}M\hat {R} (z) M^{-1}=\frac{1}{\sqrt{2}z} diag(1-i,1-i,1+i,1+i)
+\frac{z}{\sqrt{2}} diag(1+i,1+i,1-i,1-i)
\end{equation}
 The diagonal elements are complex with real trace.

\smallskip

$S14:$

\smallskip

         Here
              
\begin{equation}
\hat R =  \pmatrix{
   0 &0 &0 &q \cr
   0 &1 &0 &0 \cr
   0 &0 &1 &0 \cr
   q &0 &0 &0                                                   
              }.
\end{equation}

  The  projectors ( {\it three} even for a $4\times 4$ $\hat R$ )

\begin{equation}
P_{(0)} =  \pmatrix{
   0 &0 &0 &0 \cr
   0 &1 &0 &0 \cr
   0 &0 &1 &0 \cr
   0 &0 &0 &0                  
           }, \qquad   2P_{(\pm)} =  \pmatrix{
   1 &0 &0 &\pm1\cr
   0 &0 &0 &0 \cr
   0 &0 &0 &0 \cr
   \pm 1 &0 &0 &1                  
           }.                                         
\end{equation}

give 

\begin{equation}
\hat R = P_{(0)} + q( P_{(+)} -P_{(-)} )
\end{equation}

Baxterization gives

\begin{equation}
\hat{R}(z) = P_{(0)} + v(z)( P_{(+)} -P_{(-)} )
\end{equation}

where $v(z)$ is {\it arbitrary}. ( See $[13]$ for details.)

 One can indeed set ( with $z= e^{\theta}$, say )

\begin{equation}
v(z) = \frac{f(z)}{f(z^{-1})}, \qquad  -v(z) =  \frac{(z -z^{-1})}{(z^{-1} -z) }
\frac{f(z)}{f(z^{-1})}
\end{equation}

and factorize. But more freedom is present, as compared to $(1.1)$ and all previous
examples. Denoting  $\hat{R}(v(z))$ by $\hat{R}(v)$ one obtains

\begin{equation}
\hat {R}_{12}(v)\hat{R}_{23}(v')\hat {R}_{12}(v'')
=\hat{R}_{23}(v'')\hat{R}_{12}(v')\hat{R}_{23}(v)
\end{equation}

where $(v,v',v'')$ are mutually independent.

 Amusingly, $\hat R$ of $S14$ is diagonaized by $(4.1)$ the $\hat R$ of $S03$ giving

\begin{equation}
diag (q,1,1,-q)
\end{equation}

\section{Affine ${\cal{U}}_q (\hat{sl}_2)$:}
\setcounter{equation}{0}
\renewcommand{\theequation}{5.\arabic{equation}}

 We start below directly with the matrix $\bar R_{VV}(z)$ ( equations $(3.13)$ and
$(3.14)$ ) of $Sec.3.2$ of $[15]$. We obtain the spectral resolution and factorization (
finding back the Baxterization of $GL_q(2)$ of $Sec.2$ ). Thus, apart from a possible
overall factor, the braid matrix of ${\cal{U}}_q (\hat{sl}_2)$ is

\begin{equation}
\hat{R}(z)= P\bar R_{VV}(z) =  \pmatrix{
   1 &0 &0 &0 \cr
   0 &zc &b &0 \cr
   0 &b &c &0 \cr
   0 &0 &0 &1                                                   
              }.
\end{equation}

 where
\begin{equation}
b=\frac{(1-z)q}{(1-q^2 z)}, \qquad c=\frac{(1-q^2)}{(1-q^2 z)}
\end{equation}

Define the following basis satisfying $(1.4)$,

\begin{equation}
P_{(0)} =  \pmatrix{
   1 &0 &0 &0 \cr
   0 &0 &0 &0 \cr
   0 &0 &0 &0 \cr
   0 &0 &0 &1                  
           } ,\qquad   (q+q^{-1})P_{(\pm)} =  \pmatrix{
   0 &0 &0 &0\cr
   0 &q^{\pm1} &\pm1 &0 \cr
   0 &\pm1 &q^{\mp1} &0 \cr
   0 &0 &0 &0                  
           }.                                         
\end{equation}

   Note the specific $q$-dependence of $P_{(\pm)}$. Relabelling $(P_{(0)}+P_{(+)})$ as
$P_{(+)}$ makes the relation to $GL_q(2)$ clearer. On the other hand $(5.3)$ with $q=1$
corresponds to the basis for the $6$-vertex model $(Sec.6)$.

   Setting $q=e^h, z=e^\theta$ one can write $(5.1)$ as

  \begin{equation}
\hat R(\theta) =P_{(0)} + P_{(+)}
+\frac{sinh(h-\frac{\theta}{2})}{sinh(h+\frac{\theta}{2})}P_{(-)}                          
\end{equation}

where
\begin{equation}
\frac{sinh(h-\frac{\theta}{2})}{sinh(h+\frac{\theta}{2})}
 =\frac{(q^2-z)}{(q^2 z -1)}= \frac{z^{-\frac{1}{2}}q-
z^{\frac{1}{2}}q^{-1}}{z^{\frac{1}{2}}q- z^{-\frac{1}{2}}q^{-1}}
\end{equation}

 Factorizations of the type $(1.12)$ are now evident. Also evidently from $(5.4)$

\begin{equation}
\hat {R}(z^{-1})\hat {R}(z) =\hat {R}(-\theta)\hat {R}(\theta)= I
\end{equation}

Any supplementary overall factor, unless of the form

\begin{equation}
\frac{\rho (z)}{\rho (z^{-1})}
\end{equation}

will be incompatible with $(5.6)$. The results for $\theta \rightarrow \pm \infty$
are displayed below for comparison with the corresponding results for the $6$-vertex
$(Sec.6)$ and the $8$-vertex  $(Sec.7)$ to follow.

  For $\theta \rightarrow \infty$

\begin{equation}
\hat{R}(\theta) \rightarrow \hat R =  \pmatrix{
   1 &0 &0 &0 \cr
   0 &(1-q^{-2}) &q^{-1} &0 \cr
   0 &q^{-1} &0 &0 \cr
   0 &0 &0 &1                                                   
              }.
\end{equation}

  For $\theta \rightarrow -\infty$

\begin{equation}
\hat{R}(\theta) \rightarrow \hat R^{-1} =  \pmatrix{
   1 &0 &0 &0 \cr
   0 &0 &q &0 \cr
   0 &q &(1-q^2) &0 \cr
   0 &0 &0 &1                                                   
              }.
\end{equation}

    For $R=P\hat R$ one recognizes the the familiar lower and upper triangular $YB$
matrices of $GL_q(2)$. Thus $(5.1)$ is,indeed, a Baxterized form of $(5.8)$ and $(5.9)$
for a particular choice of basis and parametrization.

 The matrix

 \begin{equation}
M =  \pmatrix{
   1 &0 &0 &0 \cr
   0 &q &1 &0 \cr
   0 &-q^{-1} &1 &0 \cr
   0 &0 &0 &1                                                   
              }.
\end{equation}

  diagonalizes each projector giving

 \begin{equation}
M\hat R(\theta) M^{-1}= diag
(1,1,\frac{sinh(h-\frac{\theta)}{2})}{sinh(h+\frac{\theta)}{2})},1)
\end{equation}

  We have factorized the basic matrix $(5.1)$. After supplementary quasi-Hopf twists
$[15]$ one can seek again a spectral resolution to study analogous possibilities
provided that $(5.6)$ is conserved.

\section{The $6$-vertex model:}
\setcounter{equation}{0}
\renewcommand{\theequation}{6.\arabic{equation}}

  The more general $8$-vertex case is treated in $Sec.7$. But we introduce already at
this stage a basis of projectors, satisfying $(1.4)$, adequate for the $8$-vertex matrix.

 Define

\begin{equation}
2P_{1(\pm)} =  \pmatrix{
   1 &0 &0 &\pm1 \cr
   0 &0 &0 &0 \cr
   0 &0 &0 &0 \cr
   \pm1 &0 &0 &1                  
           } ,\qquad   2P_{2(\pm)} =  \pmatrix{
   0 &0 &0 &0\cr
   0 &1 &\pm1 &0 \cr
   0 &\pm1 &1 &0 \cr
   0 &0 &0 &0                  
           }.                                         
\end{equation}

 For the $6$-vertex we need only the subset

\begin{equation}
P_{(0)}=P_{1(+)} +P_{1(-)}, \qquad  P_{(\pm)} = P_{2(\pm)}
\end{equation}

  Leaving aside all well-known relations, via reparametrizations and limiting processes,
to rational affine cases, we illustrate our approach using the trigonometric
parametrization and, in particular the "ferroelectric" regime. Extensive discussions and
references can be found in the review $[16]$. ( N.B. Our $\hat R$ corresponds to $R$ in
the notation of $[16]$. See, for example, $(2.19)$ of $[16]$.)

  With our standard normalization $(Sec.2)$ in view we define (with $\gamma >0, \theta
>0$),

\begin{equation}
x=\frac {sinh \gamma}{sinh(\gamma +\theta)},\qquad  y=\frac {sinh \theta}{sinh(\gamma
+\theta)}
\end{equation}
( Though $\theta >0$ for this regime, we will consider later the limits $\theta
\rightarrow \pm \infty$. )
 The braid matrix is

\begin{equation}
\hat{R}(\theta) =  \pmatrix{
   1 &0 &0 &0 \cr
   0 &x &y &0 \cr
   0 &y &x &0 \cr
   0 &0 &0 &1                                                   
              }.
\end{equation}

 Implementing $(6.1)$ and $(6.2)$

$$\hat{R}(\theta) =  P_{(0)} + (x+y)P_{(+)}+ (x-y)P_{(-)}$$ 
\begin{equation}
  =  P_{(0)} + \frac{cosh \frac{1}{2}(\gamma -\theta)}{cosh \frac{1}{2}(\gamma
+\theta)}P_{(+)}+\frac{sinh \frac{1}{2}(\gamma -\theta)}{sinh \frac{1}{2}(\gamma
+\theta)} P_{(-)}
\end{equation}

  We have thus the structure $(1.9)$ and hence the factorization $(1.12)$. From $(6.3)$
as

\begin{equation}
\theta \rightarrow \pm \infty, \qquad  x \rightarrow 0, y\rightarrow e^{\mp \gamma}
\end{equation}

Hence the corresponding limits of $\hat{R}(\theta)$ give respectively ( see the
discussion starting with $(1.23)$ ) for the non-Baxterized $YB$ matrix

\begin{equation}
R^{\pm1}= ( P\hat R)^{\pm1} = diag (1,e^{\mp \gamma},e^{\mp \gamma},1)
\end{equation}

 This is a special class of even the simplest and the first solution $H_{3.1}$
 in the classification of $4\times 4$ $YB$ matrices $[12]$, namely

\begin{equation}
R= diag (p,q,r,s)
\end{equation}

  In view of the $8$-vertex case to follow it is convenient to choose the diagonalizer
$(M=M^{-1})$ as

 \begin{equation}
\sqrt{2} M =  \pmatrix{
   1 &0 &0 &1 \cr
   0 &1 &1 &0 \cr
   0 &1 &-1 &0 \cr
   1 &0 &0 &-1                                                   
              }.
\end{equation}

Now 
\begin{equation}
M\hat{R}(\theta) M^{-1} = diag (1,x+y,x-y,1)  
\end{equation}

 The crucial difference, for the parametrizations adopted, between $(5.1)$ and $(6.4)$
is the extra factor $z$ in the second diagonal element of $(5.1)$. This leads to the
$q$-dependent projectors in $(5.3)$ as cmpared to the elements $\pm \frac {1}{2}$ only
in $(6.1)$. Finally one is led  ( for $\theta \rightarrow \pm \infty$) to triangular 
and to diagonal $YB$ matrices    in $Sec.5$ and $Sec.6$ respectively.

\section{The $8$-vertex braid matrix:}
\setcounter{equation}{0}
\renewcommand{\theequation}{7.\arabic{equation}}

The braid matrix of the quantum affine algebra ${\cal A}_{q,p}(\hat{sl}_2)$ corresponds
to the $8$-vertex model. Some references relatively directly relevant to our purpose are
$[15],[17],[18]$ and $[19]$. These cite other basic sources. 

 Here the symmetrical structure of $(6.4)$ is generalized (with $z=e^\theta$) to

\begin{equation}
\hat{R}(z) =  \pmatrix{
   a(z) &0 &0 &d(z) \cr
   0 &c(z) &b(z) &0 \cr
   0 &b(z) &c(z) &0 \cr
   d(z) &0 &0 &a(z)                                                   
              }.
\end{equation}

 There are well-known expressions for the elements in different equivalent forms in
terms of elliptic functions. The role of a specific class of overall factors for
specific realizations of $(a,b,c,d)$ will ber commented upon later on. Implementing now
the full set $(6.1)$  with $a(z)=a$ and so on
\begin{equation}
 \hat R (z) = (a+d)P_{1(+)}+ (a-d)P_{1(-)} +(c+b)P_{2(+)}+ (c-b)P_{2(-)}
\end{equation}

 We have thus a spectral decomposition in terms of the simple basis $(6.1)$ with {\it
constant} coefficients $\pm \frac{1}{2}$. The next steps ( in order to implement $(1.9)$
and hence $(1.12)$ ) consists in explicit constructions of functions $f_{1(\pm)}(z)$
and $f_{2(\pm)}(z)$ such that $(7.2)$ satisfies $(1.1)$ for
\begin{equation}
(a\pm d)= \frac {f_{1 (\pm)} (z)}{f_{1 (\pm)} (z^{-1})}, \qquad (c\pm b)= \frac {f_{2
(\pm)} (z)}{f_{2 (\pm)} (z^{-1})}
\end{equation}

  These solutions are {\it directly} obtained from equations $(3.28)$ and $(3.29)$ of
$[15]$ in terms of infinite products

\begin{equation}
(x;a)_{\infty} = \prod_{n \geq 0} ( 1- xa^n)
\end{equation}

Noting that 

\begin{equation}
q \frac {1\pm q^{-1}z}{1\pm qz} = \pm \frac {q^{\frac{1}{2}} z^{-\frac{1}{2}} \pm
q^{-\frac{1}{2}} z^{\frac{1}{2}}}{q^{\frac{1}{2}} z^{\frac{1}{2}} \pm
q^{-\frac{1}{2}} z^{-\frac{1}{2}}}
\end{equation}

one obtains from the results cited above  ( writing $z$ for $\zeta$ and slightly
reordering the factors )

\begin{equation}
a\pm d = \frac {(\mp p^{\frac{1}{2}}q^{-1}z;p)_{\infty}(\mp
p^{\frac{1}{2}}qz^{-1};p)_{\infty}}{(\mp p^{\frac{1}{2}}q^{-1}z^{-1};p)_{\infty}(\mp
p^{\frac{1}{2}}qz;p)_{\infty}}
\end{equation}

\begin{equation}
c\pm b =\frac {(q^{\frac{1}{2}} z^{-\frac{1}{2}} \pm
q^{-\frac{1}{2}} z^{\frac{1}{2}})}{(q^{\frac{1}{2}} z^{\frac{1}{2}} \pm
q^{-\frac{1}{2}} z^{-\frac{1}{2}})} \frac {(\mp pq^{-1}z;p)_{\infty}(\mp
pqz^{-1};p)_{\infty}}{(\mp pq^{-1}z^{-1};p)_{\infty}(\mp
pqz;p)_{\infty}}
\end{equation}

 Our objectives are attained. We have arrived at $(1.9)$ and $(1.12)$.  In view of the
factored structures of $(7.6)$ and $(7.7)$ the comments $(Sec.1)$ concerning varied
possibilities in selecting $f(\theta)$ are now particularly relevant. Several factors
lead to more alternatives.

  As $\theta = ln z \rightarrow \pm \infty $, $z \rightarrow \infty$ and $z \rightarrow
0$ respectively. The extra factor in $(c\pm b)$ contributes

\begin{equation}
\biggl( \frac {q^{\frac{1}{2}} z^{-\frac{1}{2}} \pm
q^{-\frac{1}{2}} z^{\frac{1}{2}}}{q^{\frac{1}{2}} z^{\frac{1}{2}} \pm
q^{-\frac{1}{2}} z^{-\frac{1}{2}}}\biggr)_{z \rightarrow  \pm\infty} = \pm q^{\mp1}
\end{equation}

 The ratios of the infinite products ( considering the leading term in $(7.4)$ for $n
\leq k $) give  ( for $\theta \rightarrow \pm\infty $ respectively ) for both $(7.6)$
and $(7.7)$ a factor

\begin{equation}
\lim_{k \rightarrow \infty} q^{\mp2k} 
\end{equation}

 However, we have not yet implemented our standard normalization , namely
obtaining $1$ for the element $(11)$ at top left $(Sec.2)$. This is achieved by a
normalizing factor $a^{-1}$. Such a factor absorbs the limiting factor $(7.9)$. From
$(7.8)$ and $(7.9)$, for

\begin{equation}
(a',b',c',d') = a^{-1}(a,b,c,d) 
\end{equation}

\begin{equation}
lim_{\theta \rightarrow \pm \infty}(a',b',c',d') = (1, q^{\mp 1},0,0)
\end{equation}

 Note that since $a(z)$ i,e, $a(\theta)$ is itself of the form
$x(\theta)(x(-\theta))^{-1}$ such a normalization conserves the unitarity $(1.14)$
already satisfied by $(a,b,c,d)$.

   However, if one prefers to maintain the simpler symmetry of the parametrization
$(7.6),(7.7)$ one may choose to absorb the factors $(7.9)$ by a normalizing factor, say

\begin{equation}
\frac {(q^2 z;1)_\infty}{(q^2 z^{-1};1)_\infty}
\end{equation} 

This conserves $(1.14)$ and {\it again} gives the right hand side of $(7.11)$ as limits
(for $(a,b,c,d)$ normalized by $(7.12)$). This is a particularly simple choice albeit,
evidently, not unique. We do not propose to examine here normalizations adopted in the
cited sources.

   After such a normalization one obtains

\begin{equation}
\lim_{\theta \rightarrow \pm \infty}(P\hat R(\theta)) = diag
(1,q{^{\mp1}},q{^{\mp1}},1)
\end{equation}

 Thus indeed one finds again a diagonal $YB$ matrix. Compare $(6.7)$ and the comments
 preceding $(5.8)$. 

 The diagonalizer $M$ of $(6.9)$ gives now with any normalization factor $N$ and $(7.2)$

\begin{equation}
M \hat R(\theta)M^{-1} = N diag ( a+d, c+b, c-b, a-d)
\end{equation} 

\section{A nested sequence of projectors for higher dimensions:}
\setcounter{equation}{0}
\renewcommand{\theequation}{8.\arabic{equation}}

 The 8-vertex matrix has complex features due to the presence of four functions
$(a,b,c,d)$ and their realizations in terms of elliptic functions. On the other hand its
symmetry permits a spectral resolution on a basis of particularly simple symmetrical
projectors with {\it constant} elements $( \approx \pm 1)$. For $N^2 \times N^2$
matrices with $N>2$ one can construct different types of generalization of such a basis
with constant elements. One example can be easily extracted from the multistate model
presented in $Sec.4$ of $[16]$ ( where original sources are cited ). Let us consider
the simplest such case $(N=3)$.

    Let $E_{ij}$ be the matrices defined below $(1.3)$. A set of projectors satisfying
$(1.4)$ and suitable for the specral decomposition of a particular class of $9\times 9$
$\hat R (\theta)$ is

 $$ 2P_{1(\pm)} = (E_{11} +E_{99}\pm E_{19} \pm E_{91}), \qquad 2P_{2(\pm)} = (E_{22}
+E_{44}\pm E_{24} \pm E_{42}) $$
\begin{equation}
2P_{3(\pm)} = (E_{33} +E_{77}\pm E_{37} \pm E_{73}), \qquad 2P_{4(\pm)} = (E_{66}
+E_{88}\pm E_{68} \pm E_{86}) ,\qquad P_{55} =E_{55}
\end{equation}
 Generalizations for $N>3$ are not difficult to write down.

 Additional simplifications arise in $(4.1)$ of $[16]$ since the functions of $(\gamma,
\theta)$ implemented are of the $6$-vertex type. One feature should be noted. The
action of $P$ for $N=3$ interchanges the rows $(2,4),(3,7)$ and $(6,8)$. Hence when the
basis $(8.1)$ is implemented $\hat R$ and $R (=P\hat R)$ have fairly analogous structures
( as for the $4\times 4$ $6$-and $8$- vertex matrices ). Such a feature, though worth
noting, is not essential and is indeed not present in the standard cases of
$Sec.2$. 
 Instead of presenting full details concerning the above-mentioned possibility, we
briefly present another one. For $N=2$ this coincides with $(6.1)$. This basis does {\it
not} have (for $N>2$) the simple property of $(8.1)$ and its generalizations for $N>3$
under the action of $P$. But it exhibits a particularly simple canonical nested
structure. The prescription for diagonalization is also particularly simple.

         For $n=N^2 = 2l$ define

\begin{equation}
2P_{i(\pm)} = (E_{ii} +E_{n-i+1,n-i+1}\pm E_{i,n-i+1} \pm E_{n-i+1,i}),  \qquad
(i=1,2,...,l) 
\end{equation}

 For $n=2l+1$ one has in addition

\begin{equation}
P_{l+1} = E_{l+1,l+1}
\end{equation}

To diagonalize this set satisfying $(1.4)$ now define

\begin{equation}
{\sqrt 2} M ={\sqrt 2} M^{-1} = \sum_{i=1}^{l}(E_{ii} +E_{i,n-i+1}+ E_{n-i+1,i} -
E_{n-i+1,n-i+1}) + E_{l+1,l+1} 
\end{equation}

 For $n=2l$, the last term is absent.

 One obtains
\begin{equation}
MP_{i(+)}M^{-1} = E_{i,i}, \qquad MP_{i(-)}M^{-1} = E_{n-i+1,n-i+1}, \qquad (i=1,2,..l)
\end{equation}
When it is present, $P_{l+1}$ is already diagonal and commutes with $M$. Hence, if (with
the last term present only for odd $n$) and with $\varepsilon =\pm$,

\begin{equation}
\hat R(z) = \sum_{i=1}^{l} \sum_{\varepsilon} \biggl
(\frac{f_{i(\varepsilon)}(z)}{f_{i(\varepsilon)}(z^{-1})}P_{i(\varepsilon)} \biggr)
+\frac{f_{l+1}(z)}{f_{l+1}(z^{-1})}P_{l+1}
\end{equation}

 \begin{equation}
M\hat R(z) M^{-1} = \sum_{i=1}^{l}  \biggl
(\frac{f_{i(+)}(z)}{f_{i(+)}(z^{-1})}E_{ii}+\frac{f_{i(-)}(z)}{f_{i(-)}(z^{-1})}E_{n-i+1,n-i+1}
\biggr) +\frac{f_{l+1}(z)}{f_{l+1}(z^{-1})}E_{l+1,l+1}
\end{equation}

  The two sets, ( the one generalizing $(8.1)$ for all $n$ and the one given by
$(8.2),(8.3)$ ) can be shown to be related through a similarity transformation. But the
matrix of conjugation does not possess a tensored structure $G \otimes G$  ( and hence
the tensored components of the base space are not transformed individually ). The
question of existence and construction of solutions of
$(1.1)$ for the papametrization
$(8.6)$ is beyond the scope of this paper.

       Lat ue conclude with a closer look at the simplest  nontrivial case. For $N=3$,

  With a maximum number of functions in the coefficients and suppressing arguments
($x(\theta)= x$ and so on )
\begin{equation}
\hat R(\theta) =
\sum_{\varepsilon}\biggl((x+{\varepsilon} y)P_{1(\varepsilon)}+(u+{\varepsilon}
v)P_{2(\varepsilon)}
+(a+{\varepsilon}d)P_{3(\varepsilon)}+(c+{\varepsilon}b)P_{4(\varepsilon)}\biggr)
+wP_{5} 
\end{equation} 

The diagonalizer is

\begin{equation}
\sqrt 2 M = \sqrt 2 M^{-1} =  \pmatrix{
   1 &0 &0 &0&0 &0 &0&0&1 \cr
   0 &1 &0 &0&0 &0 &0&1&0 \cr                                              
   0 &0 &1 &0&0 &0 &1&0&0 \cr
   0 &0 &0 &1&0 &1 &0&0&0 \cr
   0 &0 &0 &0&\sqrt 2 &0 &0&0&0 \cr
   0 &0 &0 &1&0 &-1 &0&0&0 \cr
   0 &0 &1 &0&0 &0 &-1&0&0 \cr
   0 &1 &0 &0&0 &0 &0&-1&0 \cr
   1 &0 &0 &0&0 &0 &0&0&-1 \cr         }.
\end{equation}

  This generalizes $(6.9)$ as one moves up from the $4 \times 4$ to the $9\times 9$ case
in our sequence. Note the central element ( at (55) ) appearing for odd $n (=9)$.

\begin{equation}
M \hat R( (\theta) M^{-1} =diag( (x+y),(u+v),(a+d),(c+b),w,(c-b),(a-d),(u-v),(x-y))
\end{equation}

 If ,say, $y=0$ one can redefine $(P_{1(+)} +P_{1(-)})$ as a single projector $P_{(1)}$
and so on ( continuing to satisfy $(1.4)$ ). The number of functions available increases
with the number of projectors, but so does the number of constraints due to the braid
equation $(1.1)$. Let us note however, that for $SO_q(N)$ and $Sp_q(N)$ after fixing the
normalization one has only two functions to satisfy four complicated functional
equations $[1]$. Yet one emerges with ${\it three}$ independent solutions in both the
cases $(Secs.2,3)$. A close study of particular cases in the present context might also
lead to interesting possibilities. One recognizes $(7.2)$ and $(7.14)$ as subcases of
$(8.8)$ and $(8.10)$ respectively.

\section{Diagonalization and factorization:}
\setcounter{equation}{0}
\renewcommand{\theequation}{9.\arabic{equation}}

  Diagonalization of braid matrices was studied in $[1]$. It was used to elucidate
certain aspects of associated noncommutative spaces. Here it will be studied in the
context of factorization.

    For the $4 \times 4$ matrices ($Secs. 4,5,6,7$) the diagonalizer $M$ has been
presented for each case explicitly. In $Sec.8$ $M$ has been obtained for the nested
sequence explicitly for arbitray dimensions. The results for the lower dimensional
cases of $(A,B,C,D)_q$type algebras are collected in $App.B$. We will see how a
striking structure emerges from diagonalization of each factor in $(1.12)$. But to
start with it is worthwhile to recapitulate some basic features noted in $[1]$.

\smallskip
 $\bullet$ From $(1.6)$ and $(1.7)$ it is evident that if there exists an invertible
matrix $M$ diagonalizing $\hat R(\theta)$ it must diagonalize each projector $P_i$ ${\it
 separately}$.
\smallskip

$\bullet$  A projector, when diagonalized, can have only $+1$ or $0$ as diagonal
elements.
\smallskip

  $\bullet$ The number of unit elements on the on the diagonal is equal to the the trace
of the projector, obligatorily a positive integer. For different $P_i$ in $(1.4)$ these
elemments can never coincide due to orthogonality. In standard notations $[10]$ one has 
 \smallskip

 for $GL_q(N)$

\begin{equation}
  P_{(+)} +P_{(-)} = I_{N^2 \times N^2}           
\end{equation}
with 
\begin{equation}
2 Tr P_{(\pm)} = N(N\pm 1)            
\end{equation}

 For $SO_q(N)$ and $Sp_q(N=2n)$ one has 
\begin{equation}
  P_{(+)} +P_{(-)} +P_0 = I_{N^2 \times N^2}           
\end{equation}
 and with $\epsilon  = \pm 1$ respectively  ( as below $(3.1)$ )

\begin{equation}
2 Tr P_{(\pm)} = N(N\pm 1) \mp (\epsilon \pm 1), \qquad Tr P_0 =1           
\end{equation}
 \smallskip 
$\bullet$  If necessary, implementing a simple supplementary conjugation the elements on
the diagonal can be reordered. Exploiting this possibility we introduce the following
conventions:

 \smallskip 
 
  For $GL_q(N)$ the unit elements of $P_{(-)}$ are grouped at the top, followed by those
of $P_{(+)}$. Thus for $GL_q(3)$ and 
\begin{equation}
 \hat R(\theta)= P_{(+)} +vP_{(-)}          
\end{equation}

  once $M$ is constructed we obtain
\begin{equation}
M \hat R(\theta) M^{-1}= diag (v,v,v,1,1,1,1,1,1,)          
\end{equation}

For $SO_q$ and $Sp_q$ the chosen ordering is $ (P_0,P_{(-)} ,P_{(+)})$. Thus for
$SO_q(3)$ and 
 \begin{equation}
 \hat R(\theta)= P_{(+)} +vP_{(-)} + wP_{(0)}          
\end{equation}

\begin{equation}
M \hat R(\theta) M^{-1}= diag (w,v,v,v,1,1,1,1,1)          
\end{equation}

 Generalizations are evident.

\smallskip

$\bullet$  Since each diagonalzed $P_i$ (denoted below by $D_i$) is thus ${\it
completely}$  fixed beforehand one can ( assuming the invertibility of $M$ to be
confirmed a posteriori )  write separately for each $P_i$ with the same $M$,
\begin{equation}
 MP_i = D_i M         
\end{equation}

   Here  both $P_i$ and $D_i$ are known giving explicit ${\it linear}$
constraints on the elements of $M$. One avoids the construction of $M^{-1}$ to start
with.

\smallskip 

$\bullet$  The block structures in $(9.6)$ and $(9.8)$ and their evident generalizations
reveal the extent to which $M$ is arbitrary:

    Let $M_i$ denote a matrix of dimension $(Tr P_i \times Tr P_ i)$, with a nonzero
determinant but with otherwise ${\it arbitrary}$ elements. Then, in obvious notations,
a supplementary conjugation of $(9.6)$ by a block-diagonal $(bd)$ matrix 
\begin{equation}
(M_{(-)},M_{(+)})_{(bd)}      
\end{equation}

 and one of $(9.8)$ by

\begin{equation}
(M_{(0)},M_{(-)},M_{(+)})_{(bd)}      
\end{equation}

 leaves the diagonal forms invariant.
\smallskip

$\bullet$  The arbitrariness thus exhibited, instead of being a source of
embarrassment, provides a wide margin of maneuvre exploitable to select an $M$ with
particularly attractive properties.  We choose the following canonical feature: 
\begin {center}
 {\it mutual orthogonality of the rows of M}  . 
\end {center}
( Except for the complex, unitary $M$ of $(4.9)$ for the exotic $S03$ such an
orthogonality holds for all the cases we study.)
 
\smallskip

  Agreeable consequences are

\smallskip

  $(1)$:  The inverse of $M$ is obtained effortlessly. The prescription is:
    Take the transpose $M^T$ of $M$. Normalize each element of the column
 $j$ of $M^T$ by the {\it same} factor $c_j$ such that for each $j$

\begin{equation}
\biggl( \sum _{i} M_{ij}^2 \biggr )c_j = 1
\end{equation}

 Thus one obtains $M^{-1}$. Examples can be found in $[1]$.

\smallskip 

  $(2)$: Each row of $M$, transposed to a column provides an eigenvector of $M$ and all
         together a complete set. 

\smallskip

{\bf Consequences for factorization }:

\smallskip

     For 
\begin{equation}
 \hat{R}(\theta) =\sum_{i}\frac{f_i(\theta)}{f_i(-\theta)}P_i  
\end{equation}    

 \begin{equation}
 M\hat{R}(\theta) M^{-1} =diag
\biggl(\frac{f_1(\theta)}{f_1(-\theta)},...;\frac{f_2(\theta)}{f_2(-\theta)},...;...\biggr) 
\end{equation}

 Here the multiplicity of $f_i(\theta)$ is equal to $Tr P_i$. Note that $M$ is {\it
independent } of $\theta$. It diagonalizes each $P_i$ ( independent of $\theta$ ) and   
hence also $\hat R$.
  
  Define 

\begin{equation}
D(\theta) = diag ( f_1(\theta),...;f_2(\theta),...;...),  \qquad M(\theta) = D
(\theta)M 
\end{equation}

 Now , starting with $(1.13)$,
\begin{equation}
 \hat{R}(\theta) =(F(-\theta))^{-1}F(\theta) =
(M^{-1}D(-\theta)M)^{-1}(M^{-1}D(\theta)M) = (M(-\theta))^{-1}M(\theta) 
\end{equation}

 In each factor all $\theta$-dependence is thus again factorized in a diagonal
matrix $D(\theta)$. Some consequences will be studied in the following sections.

\section{L-operators:}
\setcounter{equation}{0}
\renewcommand{\theequation}{10.\arabic{equation}}

 Here we indicate the general features that arise as one implements our formalism in the
construction of $L$-operators. It is well-known that the FRT definitions  $[10]$ ( with
their $R^{(+)}=(PRP)$  and with $L_{2}^{\varepsilon} =PL_{1}^{\varepsilon}P$) 

\begin{equation}
(PRP) L_{1}^{\pm} L_{2}^{\pm} = L_{2}^{\pm} L_{1}^{\pm}(PRP), \qquad 
(PRP)L_{1}^{+}  L_{2}^{-}= L_{2}^{-}L_{1}^{+}(PRP)
\end{equation}

give in terms of $\hat R =PR$,

\begin{equation}
\hat R L_{2}^{\pm} L_{1}^{\pm} = L_{2}^{\pm} L_{1}^{\pm}\hat R, \qquad 
\hat RL_{2}^{+}  L_{1}^{-}= L_{2}^{-}L_{1}^{+}\hat R
\end{equation}

 Taking one more step we define ( with $\varepsilon = \pm $ below )

\begin{equation}
L_{2}^{\varepsilon} P = PL_{1}^{\varepsilon} \equiv \hat L_{\varepsilon}
\end{equation}

 when 

 $$ L_{2}^{\varepsilon} L_{1}^{{\varepsilon}'}= L_{2}^{\varepsilon} P
PL_{1}^{{\varepsilon}'} =
\hat L_{\varepsilon}\hat L_{{\varepsilon}'}$$

 and 
\begin{equation}
\hat R \hat {L}_{\varepsilon}\hat {L}_{\varepsilon}=
\hat{L}_{\varepsilon}\hat{L}_{\varepsilon}\hat R,
\qquad \hat R\hat{L}_{+}\hat {L}_{-}= \hat{L}_{-}\hat {L}_{+}\hat R
\end{equation}

   All this is before Baxterization. When the spectral parameter is introduced a more
general formulation is 

\begin{equation}
\hat R (\theta -{\theta}') \hat {L}_{\varepsilon}(\theta)\hat
{L}_{\varepsilon}({\theta}')=
\hat{L}_{\varepsilon}({\theta}')\hat{L}_{\varepsilon}(\theta)\hat R (\theta -{\theta}'),
\quad \hat R(\theta -{\theta}')\hat{L}_{+}(\theta)\hat {L}_{-}({\theta}')=
\hat{L}_{-}({\theta}')\hat {L}_{+}(\theta)\hat R(\theta -{\theta}')
\end{equation}

( For affine cases extra factors $q^{\pm c}$ can appear in the argument of $\hat R$  in
the last equation. But the above formulation suffices to illustrate our approach.)

  One can introduce a development such as 

\begin{equation}
\hat L_{\varepsilon} (\theta) =\frac{1}{\rho (\theta)} \sum _{n\geq 0}\biggl(\hat
L_{(\varepsilon,n)}e^{n\theta} +
\hat L_{(\varepsilon, -n)}e^{-n \theta}\biggr)
\end{equation}

 But for clarity in our illustrative approach let us concentrate on a particularly
simple case $[20]$. When the spectral basis on the right of $(1.9)$ has only two
projectors, $\hat R(\theta)$ can be expressed quite simply in terms of $\hat R^{\pm 1}$.
( A more general result is obtained $(3.49)$ of $[1]$.) Thus for $GL_q (N)$ from $(2.1)$
one obtains

\begin{equation}
\hat R(\theta) = \frac{e^{h+\theta} \hat R - e^{-h-\theta} \hat R^{-1}}{e^{h+\theta}-
e^{-h-\theta}}
\end{equation}

 From $(3.3)$  ( redefining $(P_{+}+  P_{-})$ as $P_1$, say ) one obtains for this new 
class of solutions
\begin{equation}
\hat R(\theta) = P_{1}  + \frac {sinh(\eta -\theta)}{sinh(\eta +\theta)}P_{0}=
\frac{e^{(\eta +\theta)} \hat R - e^{(-\eta -\theta)} \hat R^{-1}}{e^{(\eta +\theta)} -
e^{(-\eta -\theta)}}
\end{equation}

Here $\eta$ is defined for $SO_q(N)$ and $Sp_q(N)$  as in $(3.1),(3.2)$. ( Setting $\hat
R(0) =I$ one obtains the linear relation btween $\hat R$ and $\hat R^{-1}$.)  For such
cases, defining ( analogously to $(3.5.9)$ of $[20]$ , but in terms of our $\hat L$ )

\begin{equation}
\hat L(\theta) \equiv  ( e^\theta \hat L_+ -e^{-\theta}L_- ) 
\end{equation}

all the three relations $(10.4)$ can be encapsulated in the {\it single} one

\begin{equation}
\hat R(\theta -{\theta'})\hat L(\theta) \hat L({\theta}') =\hat L({\theta}') \hat
L(\theta) \hat R(\theta -{\theta'})
\end{equation}

 ( As $(10.10)$ is developed, inserting $(10.8)$ and $(10.9)$, the terms
$\hat L_{\pm}\hat L_{\mp}$ appear in "wrong order", $\hat R\hat L_{-} \hat L_{+}$ and so
on. Now, expressing $\hat R$ in trems of $\hat R^{-1}$ and vice versa one can extract
the relations $(10.4)$ with different factors depending on arbitrary $(\theta,
{\theta}')$. )

  We will use this compct forulation adapted to our special class of braid matrices
$(Sec.3)$ to illustrate the consequences of our formalism. For more general cases (
see $(10.5),(10.6)$) the basic features will be analogous along with more elaborate
sets of equations. Some indications will be given of such generalizations. Let us
however come back to our special case:

 Implementing $(10.8)$ in $(10.10)$ one obtains

$$P_0 \biggl(\hat L(\theta) \hat L({\theta}') - \hat L({\theta}') \hat L(\theta)\biggr)
P_0=0$$  

$$P_1 \biggl(\hat L(\theta) \hat L({\theta}') - \hat L({\theta}') \hat L(\theta)\biggr)
P_1=0$$  

$$P_0 \biggl((e^{\eta -\theta +{\theta}'} -e^{-\eta +\theta -{\theta}'})\hat
L(\theta) \hat L({\theta}') - (e^{\eta +\theta -{\theta}'} -e^{-\eta -\theta +
{\theta}'})\hat L({\theta}')
\hat L(\theta)\biggr) P_1=0$$  

\begin{equation}
P_1 \biggl((e^{\eta +\theta -{\theta}'} -e^{-\eta -\theta +{\theta}'})\hat
L(\theta) \hat L({\theta}') - (e^{\eta -\theta +{\theta}'} -e^{-\eta +\theta -
{\theta}'})\hat L({\theta}')
\hat L(\theta)\biggr) P_0=0
\end{equation} 

 Here , as in the general case $(10.21)$ below , the constraints are exhaustive due to
the resolution of the identity provided by $\sum P_i = I$. This aspect is evident in
the equivalent form obtained below ($(10.17),(10.18),(10.19)$) via diagonalization.

  Here the $P_i$ do not depend on $\theta$ but only on $q$. So now implementing $(10.9)$
dependence on $(\theta, {\theta}')$ becomes entirely explicit. The coefficients of
$e^{(n\theta +n' {\theta}')}$ for different $(n,n')$ must vanish separately. Only the
factors $e^{\pm \eta}$, given by $(3.1),(3.2)$ as 

\begin{equation}
tanh \eta = \sqrt {( 1 - 4 ([N-\epsilon] +\epsilon)^{-2} )}
\end{equation}

and $q$-dependent through $[N -\epsilon]$ characterizes the $L$-algebra for this specific
class of solutions.

 As emphasized in $Sec.3$, this class remains nontrivial even for $q=1$. Now for both
cases $( \epsilon = \pm 1)$, denoting $ \eta$ as $\hat {\eta}$ for $q=1$,

 \begin{equation}
tanh \hat{\eta} = \pm N^{-1}\sqrt{N^2 -4}
\end{equation}

 But the projectors $(\hat{P}_0, \hat{P}_1)$ are still different for the two cases (
$SO_q,Sp_q$ ).

 We now present the consequences of diagonalization ($Sec.9$). Both for $SO_Q (N)$ and
$Sp_q(N)$ ( remembering that $P_1 = P_+ +P_-$ ) one obtains 

 \begin{equation}
M P_0 M^{-1}= diag (1,0,....,0), \qquad M P_1 M^{-1}= diag (0,1,....,1)
\end{equation}
with
\begin{equation}
Tr P_0 =1, \qquad Tr P_1 =N^2 -1
\end{equation}

( Explicit expression for $M$ are given, in $App.B$,  only for $SO_q(3), SO_q(4)$
and $Sp_q (4)$ .)

   Define

\begin{equation}
K(\theta) = M\hat L(\theta) M^{-1}= e^{\theta}(M\hat L_{+}M^{-1}) -   e^{-\theta}(M\hat
L_{-}M^{-1}) 
\end{equation}

 Here $M$ is different for  $SO_Q (N)$ and $Sp_q(N)$. See $Sec.9$ and our
particularly simple prescription for $M^{-1}$ when the rows of $M$ are mutually
orthogonal. Our diagonalization leads to , 
    
\begin{equation}
\bigl ( K(\theta)K({\theta}') -K({\theta}') K(\theta)\bigr)_{ij} =0; \qquad (i,j) =(1,1)
,(i>1, j>1)
\end{equation}

 and for $j>1$ to

\begin{equation}
 \biggl((e^{\eta -\theta +{\theta}'} -e^{-\eta +\theta -{\theta}'})
K(\theta)  K({\theta}') - (e^{\eta +\theta -{\theta}'} -e^{-\eta -\theta +
{\theta}'}) K({\theta}')
 K(\theta)\biggr)_{1j}=0 
\end{equation} 
\begin{equation}
 \biggl((e^{\eta +\theta -{\theta}'} -e^{-\eta -\theta +{\theta}'})
K(\theta)  K({\theta}') - (e^{\eta -\theta +{\theta}'} -e^{-\eta +\theta -
{\theta}'}) K({\theta}') K(\theta)\biggr)_{j1} =0
\end{equation}

 This is the most compact form of the constraints on the $L$- operatorrs. Those on the
elements of $\hat L(\theta)$ are now obtained from

$$ \hat L(\theta) = M^{-1}\hat K(\theta) M $$

Then one can implement Gauss decomposition, if so desired, for the elements of $L^{\pm}$
to obtain results more directly comparable to those for standard cases. But {\it all
information} is encapsulated in $(10.17),(10.18)$ and$(10.19)$. All $\theta$-dependence
can be extracted as exponential factors giving the final constraints as coefficients.
For $SO_q(3)$, for example, our equations furnish the $16$ ( for $j=2,...9$ ) constraints
which involve $\eta$.

   For the more general case $( (10.5),(10.6),(10.7) )$, where 

\begin{equation}
\hat{R}(\theta) =\sum_{i}^{p}\frac{f_i(\theta)}{f_i(-\theta)}P_i
\end{equation}

the set $(10.11)$ is generalized to the following $p^2$ constraints 

 \begin{equation}
P_i\biggl ( f_{i} (\theta -{\theta}')f_{j}(-\theta +{\theta}')                       
   \hat L_{\varepsilon}(\theta) \hat L_{{\varepsilon}'}({\theta}') - f_{i} (-\theta
+{\theta}')f_{j}(\theta -{\theta}')
\hat L_{{\varepsilon}'}({\theta}') \hat L_{\varepsilon}(\theta) \biggr ) P_j =0
\end{equation}
where

$$ (\varepsilon,{\varepsilon}') = ((++), (--),(+ -)), \qquad (i,j = 1,...p )$$ 

  Diagonalization and the definition

\begin{equation}
K_{\varepsilon}(\theta) = M \hat L_{\varepsilon}(\theta) M^{-1} = \frac {1}{\rho
(\theta)}\sum_{n}\biggl (  M \hat L_{{\varepsilon} ,n} M^{-1} e^{n \theta} + M \hat
L_{{\varepsilon} ,-n} M^{-1} e^{-n \theta} \biggr )
\end{equation}

 reduces $(10.21)$ to  

 \begin{equation}
  f_{i}(\theta -{\theta}')f_{j}(-\theta +{\theta}')                       
  ( K_{\varepsilon}(\theta) K_{{\varepsilon}'}({\theta}'))_{i',j'} - f_{i} (-\theta
+{\theta}')f_{j}(\theta -{\theta}')
(K_{{\varepsilon}'}({\theta}') K_{\varepsilon}(\theta))_{i',j'}   =0
\end{equation}

 Here, for a given $(i,j)$, the ranges of $(i',j')$ are fixed by $Tr P_i$, $Tr P_j$ and
the order chosen ($Sec.9, App.B )$ for the elements unity in diagonalizing the
projectors. A simple example is provided by $ (10.18)$ and $(10.19)$. If the expansion
$(10.22)$ is a finite series ( $(10.9)$ being an extreme example ) one can extract the
limits for
$\theta$ and $(\theta -{\theta}')$ $\rightarrow \pm \infty $, since the dependence on
these parameters can be made explicit as factored coefficients. But for
a  correct extraction the functions $f_i (\theta )$ have to be properly defined ( as
noted below $(1.33)$ ).

\section{Transfer matrices and diagonalization:}
\setcounter{equation}{0}
\renewcommand{\theequation}{11.\arabic{equation}}
 
\smallskip

 {\bf General formulation } :

\smallskip

      We start by introducing notations analogous to those of $Sec.10$ for the
row-to-row transfer matrix $T^{(L)} (\theta)$, satisfying 

\begin{equation}
\hat R(\theta - {\theta}')  \bigl (T^{(L)} (\theta) \otimes T^{(L)} ({\theta}') \bigr )
=\bigl (T^{(L)}  ({\theta)}' \otimes T^{(L)} (\theta) \bigr ) \hat R(\theta - {\theta}')
\end{equation}

 Here, apart from evident analogies ( since we have again a class of $L$-functions )
specific features arise concerning the component blocks of $T^{(L)}$. The dimensions of
the blocks increase with the length of the row according to standard prescriptions.

 Matrix multiplication for a $N^2 \times N^2$ matrix $\hat R (\theta)$  is defined by
labelling $T^{(L)}$ for any $L$ by $N^2$ blocks. If $I$ be the $N \times N$
 unit matrix,

 $$  T^{(L)} (\theta) \otimes T^{(L)} ({\theta}') =(T^{(L)} (\theta) \otimes I)(I
\otimes T^{(L)} ({\theta}')) = (P(I
\otimes T^{(L)} (\theta)P)(P(T^{(L)} ({\theta}') \otimes I)P)$$

\begin{equation}
= (P(I\otimes T^{(L)} (\theta))((T^{(L)} ({\theta}') \otimes I)P ) \quad = (P
T_2^{(L)}(\theta) )( T_1^{L}({\theta}' ) P ) \quad = {\hat T}^{(L)} (\theta) 
\hat {T}^{(L)} ({\theta}')
\end{equation}

 where  $$ \hat {T}^{(L)} \equiv P T_2^{(L)} = T_1^{(L)}P $$

are $N^2 \times N^2$ matrices in terms of blocks of $T^{(L)}$.

 Thus,

\begin{equation}
\hat R(\theta - {\theta}')  \bigl (\hat {T}^{(L)} (\theta) \hat {T}^{(L)} ({\theta}')
\bigr ) =\bigl (\hat{T}^{(L)}  ({\theta)}' \hat {T}^{(L)} (\theta) \bigr ) \hat R(\theta
- {\theta}')
\end{equation}

    Thus the mixture of matrix multiplication and tensor product in $(11.1)$ has been
rephrased as matrix multiplications of $\hat R$ and $\hat T$. Instead of $(T_1,T_2)$ the
same $\hat T$ now appears throughout.

   For $(1.9)$, namely

\begin{equation}
\hat{R}(\theta) =\sum_{i}^{p}\frac{f_i(\theta)}{f_i(-\theta)}P_i
\end{equation}

one obtains, as in $Sec.10$, a complete set of $p^2$ constraints
 \begin{equation}
P_i\biggl ( f_{i} (\theta -{\theta}')f_{j}(-\theta +{\theta}')                       
   \hat{T}^{(L)} (\theta) \hat{T}^{(L)}({\theta}') - f_{i} (-\theta
+{\theta}')f_{j}(\theta -{\theta}')
\hat {T}^{(L)}({\theta}') \hat {T}^{(L)}(\theta) \biggr ) P_j =0
\end{equation}

 Using the diagonalizer $M$ of the braid matrix ( $Sec.9, App.B$ )  define 

\begin{equation}
\hat{R}_{d}(\theta)= M\hat{R}(\theta)M^{-1}, \qquad \hat{K}^{(L)}(\theta)=
M\hat{T}^{(L)}(\theta)M ^{-1}
\end{equation}

 one obtains from $(11.3)$, in terms of the {\it diagonal} matrix $\hat R_{d}$

\begin{equation}
\hat R_{d}(\theta - {\theta}')  \bigl (\hat {K}^{(L)} (\theta) \hat {K}^{(L)} ({\theta}')
\bigr ) =\bigl (\hat{K}^{(L)}  ({\theta)}' \hat {K}^{(L)} (\theta) \bigr ) \hat
R_{d}(\theta - {\theta}')
\end{equation}

 This corresponds to

 \begin{equation}
 f_{i} (\theta -{\theta}')f_{j}(-\theta +{\theta}')                       
   \bigl(\hat K^{(L)}(\theta) \hat K^{}(L)({\theta}') \bigr )_{i'j'} - f_{i} (-\theta
+{\theta}')f_{j}(\theta -{\theta}')
\bigl (\hat K^{(L)}({\theta}') \hat K^{(L)}(\theta) \bigr )_{i'j'}   =0
\end{equation}

  We have explained  below $(10.23)$ of $Sec.10$ how the domain of $(i'j')$ depend on
the conventions adopted for the diagonalizations of the projectors $P_i$ and $P_j$. It
was also pointed out before that such a set of constraints is exhaustive. 

 The elements of $\hat K^{(L)}(\theta)$ are linear combinations of those of             
$ T^{(L)}(\theta)$, the coefficients being independent of $\theta$ ( since $M$ is so ).
The bilinear algebraic relations between  due to $(1.1)$ between $T^{(L)}_{ij}(\theta)$
attain their simplest form in $(11.8)$ in terms of these linear combinations.
Construction of a "$ \hat K$-basis" ( a complete set of states specifically adapted to
the action of the blocks of $\hat K$ ) would permit a full exploitation of $(11.8)$.

\smallskip

 {\bf Particular cases:}

\smallskip

                 We now consider two particular cases. The first one is chosen because
it is familiar and extensively studied. The content of $(11.8)$ for the $6$-vertex case
can be compared to well-known results. ( See $[7,8,9,16]$ and basic sources cited in
these references.) The second one is chosen as a relatively simple but new example of a
multistate model. It corresponds to our special class of solutions ($Sec.3$) for
$SO_q(3)$. This can be compared to a different class of multistate models $[16]$.

\smallskip

$\bullet$  The $6$-vertex case:

\smallskip

       Inserting in $(11.6)$ the results of $Sec.6$ with $M$ given by $(6.9)$, 

\begin{equation}
\hat R_{d}(\theta - {\theta}')  = diag (1,u,v,1)
\end{equation}

where

\begin{equation}
u= \frac {cosh \frac{1}{2}(\gamma - \theta + {\theta}') }{cosh \frac{1}{2}(\gamma +
\theta - {\theta}')} \qquad v= \frac {sinh \frac{1}{2}(\gamma - \theta + {\theta}')
}{sinh \frac{1}{2}(\gamma +
\theta - {\theta}')}
\end{equation}
      
 ( Other interesting choices of $M$ are possible. But $(6.9)$ is adequate for our
present purposes.) Now let 

\begin{equation}
\hat{T}^{(L)}(\theta) =  \pmatrix{
  A &B \cr
   C &D }
          .                                                        
\end{equation}

where each entry is a $2^L \times 2^L$ block obtained according to standard
prescriptions ( e.g. $Secs. 2,3$ of $[16]$ ). From $(6.9),(11.2)$ and $(11.6)$ one
obtains

\begin{equation}
2\hat{K}^{(L)}(\theta) =  \pmatrix{
  A+D & B+C & B-C & A-D \cr
  B+C & A+D & D-A & C-B \cr
  B-C & A-D & -A-D & -B-C \cr
  A-D & B-C & B+C & A+D }.
 \end {equation} 

\begin{equation}
 = (A+D) \pmatrix{
 s_0&0 \cr
   0&-s_3 } +(A-D) \pmatrix{
 0&-s_2 \cr
   s_1&0} + (B+C)\pmatrix{
 s_1&0 \cr
   0&s_2 } +(B-C)\pmatrix{
 0&s_3 \cr
   s_0&0}
          .                                                        
\end{equation}

where 

$$ s_0 = \pmatrix{
 1&0 \cr
   0&1}, \quad s_3 = \pmatrix{
 1&0 \cr
   0&-1},\quad s_1 = \pmatrix{
 0&1 \cr
   1&0}, \quad s_2 = \pmatrix{
 0&-1 \cr
   1&0} $$

 Note that 

  $$ Tr \hat K^{(L)}(\theta) = Tr (A+D) = Tr T^{(L)}(\theta) $$

and if $V$ is an eigenvector of $\hat K(\theta)$ then $M^{-1}V$ is one
of $\hat T^{(L)}(\theta)$.

   Now $(11.8)$ reduces to

\begin{equation}                      
 \biggl(\hat K^{(L)}(\theta) \hat K^{}(L)({\theta}')  - 
x\hat K^{(L)}({\theta}') \hat K^{(L)}(\theta) \biggr )_{ij}   =0
\end{equation}

where $(u,v)$ being given by $(11.10)$, one obtains $x$ as follows for values of $(i,j)$
indicated at right :

  $$ x=1 , \qquad (i,j)= (1,1),(2,2), (3,3), (4,4),(1,4), (4,1); $$
  $$x=u, \quad (i,j)= (1,2),(4,2); \qquad  x= u^{-1}, \quad (i,j) = (2,1), (2,4) ; $$
  $$ x=v, \quad (i,j) = (1,3),(4,3); \qquad x= v^{-1}, \quad (i,j)=(3,1), (3,4);$$
\begin{equation}
  x= \frac {u}{v},\quad (i,j)= (3,2); \qquad  x= \frac {v}{u},\quad (i,j)= (2,3)
\end{equation}

 Thus we obtain the simplest form of the constraints implied by $(11.1)$ or $(11.3)$.

\smallskip

       $\bullet$  A special class of multistate models ($SO_q(3)$ example):

\smallskip

   We consider now the class of braid matrices presented in $Sec.3$ and explicitly
diagonalized for $SO_q(3)$, as well as for $SO_q(4)$, in $App.B$. ( We consider only
$SO_q(N)$. For $Sp_q(N)$ certain states have negetive weights. ) The precise way in which
the model is "nonminimal", with more than two possible states per link, will be
explained at the end by comparing it with another class of models studied in $Sec.4$ of
$[16]$ ( where original sources are cited ).

 We start with $(B.8),(B.9)$ and $(B.10)$. Using the notations of $(10.8)$ the
braid matrix is ( with $\eta$ defined in $Sec.3$ )

\begin{equation}
\hat R(\theta) = P_{1}  + \frac {sinh(\eta -\theta)}{sinh(\eta +\theta)}P_{0} = I +
\biggl( \frac {sinh(\eta -\theta)}{sinh(\eta +\theta)} - 1 \biggr)P_{0}
\end{equation}

For $N=3$ (see $Sec.4$ of $[1]$),

\begin{eqnarray}
(q+q^{-1}+1) P_0 = q^{-1} E_{11}\otimes E_{33}+  q^{-\frac{1}{2}} E_{12}\otimes E_{32}+
 E_{13}\otimes E_{31} \nonumber  \\ + q^{-\frac{1}{2}} E_{21}\otimes E_{23}+
E_{22}\otimes E_{22}+ q^{\frac{1}{2}} E_{23}\otimes E_{21} \nonumber  \\  + E_{31}\otimes
E_{13}+ q^{\frac{1}{2}} E_{32}\otimes E_{12}+ q E_{33}\otimes E_{11}
\end{eqnarray}

Setting $\epsilon =1$ and $N=3$ in $(3.1)$  and $(3.2)$

\begin{equation}
2 cosh \eta = (q+q^{-1}+1)
\end{equation}

and for $M$ given by $(B.10)$

\begin{equation}
M \hat R(\theta) M^{-1} = diag \biggl(\frac {sinh(\eta -\theta)}{sinh(\eta
+\theta)},1,1,1,1,1,1,1,1 \biggr )
\end{equation}

 Hence $(11.8)$ reduces to

 \begin{equation}                      
 \biggl(\hat K^{(L)}(\theta) \hat K^{}(L)({\theta}')  - 
 K^{(L)}({\theta}') \hat K^{(L)}(\theta) \biggr )_{ij}   =0
\end{equation}

 for $(i,j)=(1,1)$ and for $(i>1,j>1)$. For $j=(2,3,...,9)$ one obtains

  \begin{equation}                      
 \biggl(sinh(\eta -\theta +{\theta}')\hat K^{(L)}(\theta) \hat K^{(L)}({\theta}')  - 
sinh(\eta +\theta -{\theta}')\hat K^{(L)}({\theta}') \hat K^{(L)}(\theta) \biggr )_{1j}  
=0
\end{equation}

\begin{equation}                      
 \biggl(sinh(\eta +\theta -{\theta}')\hat K^{(L)}(\theta) \hat K^{(L)}({\theta}')  - 
 sinh(\eta -\theta +{\theta}')\hat K^{(L)}({\theta}') \hat K^{(L)}(\theta) \biggr
)_{j1}   =0
\end{equation}
 
   Thus we obtain the complete set of constraints in the simplest and the most compact
form. The foregoing structure is directly generalizable to all $N$. But for $N>4$ one
 has either to construct  the corresponding $M$ or to use $(11.5)$ with $(11.16)$. A
study of the $\hat K$-basis adapted to the foregoing set of constraints is beyond the
scope of this paper. We conclue with some comments and comparisons.

  In $6$- or $8$- vertex models $2$ states are possible for each link. But non-zero
Boltzmann weights are associated to a subset of the $2^4$ elements of the braid matrix.
When $3$ states are possible per link ( of a plane lattice ) one can implement  a
$9\times 9$ matrix  ( with $3^4$ elements ) attributing again non-zero weights to a
subset of the possible states, corresponding to the non-zero elements of the matrix. 

    The number of non-zero elements in $(11.16)$ is $15 (=3(2.3 - 1))$. For $SO_q(N)$
this number, for our class of solution of $Sec.3$, is $N(2N-1)$. In $Sec.4$ of $[16]$ a
class of models is studied where one has precisely the same number of non-zero weights
out of $N^4$ elements ( with the symbol $q$ for our $N$ ). Apart from this feature, the
block structure in $(4.12)$ of $[16]$ ( namely, $(11),(1j),(j1),(ij)$ with $(i>1,j>1)$ )
corresponds also to the structure of our set of constraints  $( (11.20),(11.21),(11.22)
)$ for $\hat K(\theta)$. We will not study here the different possibilities concerning
block srtuctures but emphasize that in spite of the foregoing features the two classes
are basically different. For $N=2$ that of $[16]$ reduces to the $6$-vertex case,
whereras our special class does not exist. Ours is obtained for $SO_q(N)$ and and for
all $N \geq 3$ can be diagonalized to the form $(11.19)$ with $(N^2 -1)$ unit elements. 

\smallskip

  Finally let us note the situation for $q=1$. As pointed out at the end of $Sec.3$,
this class of $\hat R(\theta)$ remains nontirvial, even quite interesting, for $q=1$.
There is, of course, additional simplicity. Thus denoting $(P_0, \eta )$ for $q=1$ by 
$(\hat {P}_0, \hat{\eta} )$ one obtains for  $SO(3)$, for example,

\begin{eqnarray}
3 \hat {P}_0 =  E_{11}\otimes E_{33}+  E_{12}\otimes E_{32}+
 E_{13}\otimes E_{31} \nonumber  \\ +  E_{21}\otimes E_{23}+
E_{22}\otimes E_{22}+ E_{23}\otimes E_{21} \nonumber  \\  + E_{31}\otimes
E_{13}+ E_{32}\otimes E_{12}+  E_{33}\otimes E_{11}
\end{eqnarray}

and 

 \begin{equation}
cosh \hat {\eta} =\frac{3}{2}, \qquad sinh \hat {\eta} =\frac{\sqrt 5}{2}
\end{equation}

 It is amusing to note the relation of $ sinh \hat {\eta}$ with the Golden Mean.

\section{Remarks:}
\setcounter{equation}{0}
\renewcommand{\theequation}{12.\arabic{equation}}

The first essential step in our approach has been the spectral decomposition of the
braid matrices $\hat R(\theta)$, obtaining the coefficient of the projector $P_i$ in the
form of a ratio $\frac {f_i (\theta)}{f_i (-\theta)}$. 

 The next major step was diagonalization. But here again the projectors play a basic
role. The fact that the same matrix $M$ must diagonalize each projector appearing in the
decomposition permits a systematic exctraction of the necessary constraints on $M$ and
also the exploitation of the remaining freedom in an efficient fashion. This is explained
in detail in $Sec.9$. Along with diagonalization a remarkable new feature arises inthe
factorization. In each factor the dependence on the spectral parameter $\theta$ is again
factorized out as a diagonal matrix.

   Various directions opened up deserve further exploration. The present work, despite
its length, stops short at various points. Applications of our formalism to
$L$-operators and to transfer matrices have merely been adumbrated. Thus the
introduction of basis states specifically adapted to the form of the constraints
obtained ( "$\hat K$-basis") can be particularly helpful. For the braid matrices of
$Sec.2$ the diagonalizers have been constructed explicitly ($App.B$) only for lower
dimensions. While for $GL_q(N)$ the general prescription should not be difficult to
obtain, for $SO_q(2n +1),SO_q(2n)$ and $Sp_q(2n)$ one has, among other things, to obtain
the mutually orthogonal $(2n+1)$-plets and $2n$-plets ( $App.B $). The problem of
solving the braid equation $(1.1)$ implementing the nested sequence of projectors of
$Sec.8$ has not been addressed. 

   We hope to study elsewhere some of the aspects mentioned above.

\section {APPENDIX A :Comparison with triangular factorization}
\setcounter{equation}{0}
\renewcommand{\theequation}{A.\arabic{equation}}

  We compare here our factorization schemes with that proposed by Maillet et al.
$[7,8]$. We start with notations and general features.

  For our $R_q(\theta)$, where $q=e^{h}$, define

\begin{equation} 
z_1 = e^{(h-\theta)}, \qquad z_2 = e^{(h+\theta)}
\end{equation}

 Then as 
\begin{equation} 
\theta \rightarrow -\theta, \qquad (z_1,z_2) \rightarrow (z_2,z_1)
\end{equation}

 The unitarity $(1.14)$ is now ( in terms of $R =P\hat R$ )

\begin{equation} 
R_{21}(z_2,z_1) R_{12} (z_1,z_2) = I
\end{equation}

 In $[7,8]$ the proposed factorization is

\begin{equation} 
 R_{12} (z_1,z_2) = (F_{21}(z_2,z_1))^{-1} F_{12} (z_1,z_2) 
\end{equation}
where the aim is to obtain lower triangular $F_{12}$.
Our $(1.12)$ corresponds ( implicitly with a {\it different} $F$ )to 

\begin{equation} 
 R_{12} (z_1,z_2) = (F_{21}(z_2,z_1))^{-1}P F_{12} (z_1,z_2) 
\end{equation}

Note the presence of $P$ in $(A.5)$. In a complimentary fashion, for the braid matrix
$(A.4)$ leads to 

\begin{equation} 
\hat R_ (z_1,z_2) = (F(z_2,z_1))^{-1}P F (z_1,z_2) 
\end{equation}

as compared to our  $(1.12)$  

\begin{equation} 
\hat R_ (z_1,z_2) = (F(z_2,z_1))^{-1} F (z_1,z_2) 
\end{equation}

 This last form permits us to fully exploit the spectral decomposition.

\smallskip

 Let us note the following features:

\smallskip

 $\bullet$  Given a $R(z_1,z_2)$ one has to extract $F(z_1,z_2)$ from $(A.4)$. For
higher dimensions this ( and in particular the explicit construction of $F^{-1}$ ) is
difficult. So, assuming invertibility, the authors start from  

\begin{equation} 
(F_{21}(z_2,z_1)) R_{12} (z_1,z_2) =  F_{12} (z_1,z_2) 
\end{equation}

  For the $4 \times 4$ matrix of the $6$-vertex  model ( see our $Sec.6$ ) explicit
triangular factors in $(A.4)$ are obtained. For constructing transfer matrices
"partial" $F$-matrices are defined.

    Given our Baxterization, {\it our type of factorization is obtained effortlessly as
a byproduct}. In $[1]$ the forms $(1.9)$ (reproduced here in $Sec.2$ and $Sec.3$) were
obtained in a quest for elegance. Factorization was not a goal.

\smallskip

  $\bullet$  The limiting cases $ \theta \rightarrow  \pm \infty $ corresponds to
$(z_1,z_2)  \rightarrow  (0,\infty) ,(\infty,0) $ respectively. We have systematically
extracted the standard (non-Baxterized) braid matrices as  the corresponding limits of
the $\theta $-dependent braid matrices ($ Secs. (2,3,...,7)$) and explained in what sense
precisely ($Sec.1$) they can still be considered to be factorized.

\smallskip

 $\bullet$  It is instructive to compre different types of factorization explicitly for
the simple example of the $6$-vertex matrix. 

  Setting $\theta = (\lambda -\mu)$ in $(69)$ of $[7]$ one obtains from $(89)$ and
$(90)$ of $[7]$ , ,using a block-diagonal notation,

\begin{equation}
F_{12} (\theta)= (1,B,1)_{bd} 
\end{equation}
 where the $2\times 2$ block $B$ is 

\begin{equation}
B =  \pmatrix{
  1 &0 \cr
   \frac {sinh \eta}{ sinh (\eta + \theta)} & \frac {sinh \theta}{ sinh (\eta + \theta)}}
          .                                                        
\end{equation}

  Using the results of $Sec.6$ one obtains for our $F$,

\begin{equation}
F (\theta)= (1,C,1)_{bd} 
\end{equation}
where 
\begin{equation}
C =  \pmatrix{
  e^{\frac{1}{2}(\gamma -\theta)} & e^{-\frac{1}{2}(\gamma -\theta)} \cr
  e^{-\frac{1}{2}(\gamma -\theta)} &  e^{\frac{1}{2}(\gamma -\theta)}}
          .                                                        
\end{equation}

\smallskip

$\bullet$  Define the diagonal matrices

\begin{equation}
D(\theta)= diag ( 1, cosh \frac{1}{2}(\gamma -\theta),sinh\frac{1}{2}(\gamma
-\theta), 1) 
\end{equation}

\begin{equation}
(D(-\theta))^{-1}-= diag ( 1, (cosh \frac{1}{2}(\gamma
+\theta))^{-1},(sinh\frac{1}{2}(\gamma +\theta))^{-1}, 1) 
\end{equation}

 Now, using the $M$ of $(6.9)$, $\hat R(\theta)$ of $(6.5)$ and inverting $(6.10)$, one
can write ( with an $M$ independent of $\theta$ )

\begin{equation}
\hat R(\theta) = (M^{-1}(D(-\theta))^{-1}) (D(\theta) M) \equiv
(M(-\theta))^{-1}M(\theta)
\end{equation}

 Now in each factor $M(\theta)$ all ${\theta}$-dependence is again factorized as a
diagonal matrix. This is a general feature of this approach.

 \smallskip

  We have compared different types of factorization. One can hope to implement
fruitfully in different contexts their complementary features such as those indicated
above.

\section {APPENDIX B: Explicit diagonaliaztions}

\setcounter{equation}{0}
\renewcommand{\theequation}{B .\arabic{equation}}

  In $Sec.9$ general aspects of diagonalization of braid matrices have been presented.
Here we give explicit expressions for matrices $M$ diagonalizing $\hat R (\theta)$ for $
GL_q(2),GL_q(3),SO_q(3),SO_q(4)$ and $Sp_q(4)$.

   The result for $GL_q(2)$ effectively appears in $Sec.5$ in a form suited to the
context. Here we give an equivalent form consistent with the canonical convention of
$Sec.9$ ( see $(9.5)$ and $(9.6)$ ). Though we stop with $GL_q(3)$, one can see the
general structure of $M$ for $GL_q(n)$ emerging. For the orthogonal and the symplectic
cases the situation will be discussed at the end.

 \begin {center}    
{\it In each case below the rows of $M$ will be mutually orthogonal.}  
\end {center}

 Hence $M^{-1}$, always given by the prescription $(9.12)$, will not be displayed
explicitly.

  For $GL_q(n)$ we adopt ( with the matrices $E_{ij}$ defined below $(1.3)$ ) the
normalization

\begin {equation}
R_q = \sum _{i} E_{ii} \otimes E_{ii} + q^{-1} \sum _{i \neq j} E_{ii} \otimes E_{jj}
+ ( 1-q^{-2})\sum _{j >i} E_{ij} \otimes E_{ji}
\end{equation}
 
 The, braid matrix is
 
\begin {equation}
\hat R_q = P R_q = P_{(+)} -q^{-2} P_{(-)}
\end{equation}

 With the notations of $Sec.2$,

\begin {equation}
\hat R_q (\theta)= P_{(+)} + \frac {sinh ( h- \theta )}{ sinh (h + \theta )} P_{(-)}
\quad \equiv  P_{(+)} + v (\theta ) P_{(-)}
\end{equation}

 The projectors depend on $q(=e^h)$ but not on $\theta$. 
  
  For $n=2$ and 

\begin {equation}
M=(E_{12} -q^{-1}E_{13}) +(E_{22} +qE_{23}) +E_{31} +E_{44}
\end{equation}

\begin {equation}
M\hat R_q (\theta) M^{-1}= diag ( v(\theta),1,1,1)
\end{equation}

 For $n=3$ and 
\begin {eqnarray}  
M=(E_{12} -q^{-1}E_{14}) +(E_{23} -q^{-1}E_{27}) +(E_{36} -q^{-1}E_{38}) \nonumber \\ +
(E_{42} +qE_{44})+(E_{53} +qE_{57})+(E_{66} +qE_{68}) \nonumber \\ + E_{71}+E_{85} +
E_{99}
\end{eqnarray}
 
\begin {equation}
M\hat R_q (\theta) M^{-1}= diag ( v( \theta),v(\theta),v(\theta),1,1,1,1,1,1)
\end{equation}

   The emerging general sructure of $M$ for $GL_q(n)$ is as follows:

  There are $\frac{1}{2}n(n-1)$ rows with two nonzero elements $(1, -q^{-1})$, suitably
shifted horizontally in successive rows to assure mutual orthogonality. Then there are 
  $\frac{1}{2}n(n-1)$ rows with two nonzero elements $(1,q)$ in the corresponding
columns ( as in $(B.6)$ ). Then there are $n$ rows with a single nonzero element $1$ in
otherwise empty columns.

 \smallskip

    For $SO_q(3)$, $SO_q(4)$ and $Sp_q(4)$

\begin {equation}
\hat R_q (\theta)  = P_{(+)} +v(\theta) P_{(-)}  +w(\theta) P_{(0)}
\end{equation}

 For the orthogonal case the three possibilities for $v(\theta)$ and $w(\theta)$ are
given  (with $n=3$ and $n=4$ respectively ) by $(2.2)$,$(2.3)$ and also $(3.3)$ with
$\epsilon =1$. For the symplectic case the relevant equations are $(2.4),(2.5)$ (N.B.with
$n=2$ there ) and $(3.3)$ with $\epsilon =-1$.

   For $SO_q(3)$ define

 \begin {equation}
s= -q^{-\frac{1}{2}}(1-q), \quad t=-q^{-\frac{3}{2}}(1+q)
\end{equation}

 Now,

\begin {eqnarray}  
M=(E_{13} +q^{\frac{1}{2}}E_{15} +qE_{17}) +(E_{22} -qE_{24}) +(E_{36}
-qE_{38}) \nonumber \\
 + (E_{43} +sE_{45}-E_{47}) +E_{51}+(E_{62} +q^{-1}E_{64}) \nonumber \\
 +( E_{73}+tE_{75}+q^{-2}E_{77}) + (E_{86}+q^{-1}E_{88})+ E_{99}
\end{eqnarray}

 gives 

\begin {equation}
M\hat R_q (\theta) M^{-1}= diag (w(\theta), v( \theta),v(\theta),v(\theta),1,1,1,1,1)
\end{equation}

For $SO_q(4)$

\begin{eqnarray}
M= (E_{14}+qE_{17}+q E_{1,10}+q^2 E_{1,13})+(E_{24}+qE_{27}-q^{-1} E_{2,10}- E_{2,13})
\nonumber \\ +(E_{32}-qE_{35})+(E_{43} -qE_{49})+((E_{58}- qE_{5,14})+(E_{6,12} -
qE_{6,15}) \nonumber \\^+ (E_{74}-q^{-1}E_{77}+qE_{7,10}-E_{7,13}) +
E_{81}+(E_{92}+q^{-1}E_{95})\nonumber \\+(E_{10,3}+q^{-1}E_{10,9}) 
+ ( E_{11,4}-q^{-1}E_{11,7}-q^{-1}E_{11,10}+q^{-2}E_{11,13}) \nonumber \\+
(E_{12,8}+q^{-1}E_{12,14})
+(E_{13,12}+q^{-1}E_{13,15})+ E_{14,11} +E_{15,6}+E_{16,16}
\end{eqnarray}

gives 
\begin {equation}
M\hat R_q (\theta) M^{-1}= diag (w(\theta),
v(\theta),v(\theta),v(\theta),v(\theta),v(\theta),v(\theta),1,1,1,1,1,1,1,1,1)
\end{equation}

For $Sp_q(4)$

\begin{eqnarray}
M= (E_{14}+qE_{17}-q^3 E_{1,10}-q^4 E_{1,13})+(E_{24}-q^{-1}E_{27}+q E_{2,10}-
E_{2,13})
\nonumber \\ +(E_{32}-qE_{35})+(E_{43} -qE_{49})+((E_{58}- qE_{5,14})+(E_{6,12} -
qE_{6,15}) \nonumber \\ (E_{74}+qE_{77}+q^{-3}E_{7,10}+q^{-2}E_{7,13}) +
( E_{84} -q^{-1}E_{87} -q^{-1}E_{8,10}+q^{-2}E_{8,13})
\nonumber \\+(qE_{92}+E_{95})
+(qE_{10,3}+E_{10,9})  + (qE_{11,8}+E_{11,14}) \nonumber \\+
(qE_{12,12}+E_{12,15})+E_{13,1}+ E_{14,6} +E_{15,11}+E_{16,16}
\end{eqnarray}

 gives
 \begin {equation}
M\hat R_q (\theta) M^{-1}= diag (w(\theta),
v(\theta),v(\theta),v(\theta),v(\theta),v(\theta),1,1,1,1,1,1,1,1,1,1)
\end{equation}

 Note that the multiplicity of $v(\theta)$ is $6$ in $(B.13)$ and $5$ in $(B.15 )$.
       
 Other examples of $M$ can be found near the ends of $Secs.(4,5,6,7)$. In $Sec.8$ $M$
is obtained for for arbitrary dimensions.

  For $GL_q(n)$ one encounters as elements of different rows , apart from singlets
(unity), only the mutually orthogonal doublets 

$$ (1,-q^{-1}),\quad (1,q)$$

 But for $SO_Q(3)$ one has ( implementing $(B.9)$ ) also the mutually orthogonal triplets
 
   $$(1,q^{\frac{1}{2}},q),\quad (1,s,-1),\quad (1,t,q^{-2})$$

 For $SO_q(4)$ and $Sp_q(4)$, respectively, one similarly encounters the the mutually
orthogonal quadruplets

$$ (1,q,q,q^2),(1,q,-q^{-1},-1),(1,-q^{-1},q,-1),(1,-q^{-1},-q^{-1},q^{-2})$$

$$ (1,q,-q^3,-q^4),(1,q^{-1},q,-1),(1,q,q^{-3},q^{-2}),(1,-q^{-1},-q^{-1},q^{-2})$$

 In $Sec.7$ of $[1]$ the relations of such multiplets with particular types of
$q$-deformed surfaces ( spheres, hyperboloids) have been pointed out. In constructing
$M$ for $SO_q(N)$ and $Sp_q(N)$ for higher dimensions a key feature would be the
general srtucture   of the corresponding $N$ mutually orthogonal $N$-plets.

\bigskip

\bigskip

 This paper is dedicated to the memory of P's elegant participations in my work.

\bigskip

\bigskip

\bibliographystyle{amsplain}

\end{document}